\documentclass[final]{svjour3}

\usepackage{a4wide}
\usepackage{amssymb,amsmath}
\usepackage{graphicx}
\usepackage[english]{babel}
\usepackage{enumitem}
\usepackage{relsize}

\def\rr{\mathbb{R}} 
\def\sF{\mathcal{F}} 
\def\Pr{\mathbb{P}} 
\def\EX{\mathbb{E}}	
\def\dif{\mathrm{d}} 
\def\landau{\mathcal{O}} 
\def\part{\mathcal{T}_{\mathbf{h}}^N} 
\def\e{\textnormal{e}} 

\def\eps{\epsilon}
\def\ntilde{\widetilde}

\def\nbar{\overline}

\title{Laws of large numbers and Langevin approximations for stochastic neural field equations
}
\titlerunning{Limit theorems for stochastic neural field equations}        

\author{Evelyn Buckwar \and Martin G.~Riedler}

\institute{E.~Buckwar \and M.~G.~Riedler\at
              Johannes Kepler University \\
	      Institute for Stochastics\\
              Tel.: +43 732 2468 4160 \\
              Fax: +43 732 2468 4162\\
              \email{evelyn.buckwar@jku.at} \and \email{martin.riedler@jku.at}           
}

\date{Received: date / Accepted: date}

\numberwithin{equation}{section}
\numberwithin{theorem}{section}
\numberwithin{corollary}{section}
\numberwithin{remark}{section}

\begin{document}

\maketitle

\begin{abstract} In this study we consider limit theorems for microscopic stochastic models of neural fields. We show that the Wilson-Cowan equation can be obtained as the limit in probability on compacts for a sequence of microscopic models when the number of neuron populations distributed in space and the number of neurons per population tend to infinity. Though the latter divergence is not necessary. This result also allows to obtain limits for qualitatively different stochastic convergence concepts, e.g., convergence in the mean. Further, we present a central limit theorem for the martingale part of the microscopic models which, suitably rescaled, converges to a centered Gaussian process with independent increments. These two results provide the basis for presenting the neural field Langevin equation, a stochastic differential equation taking values in a Hilbert space, which is the infinite-dimensional analogue of the Chemical Langevin Equation in the present setting. On a technical level we apply recently developed law of large numbers and central limit theorems for piecewise deterministic processes taking values in Hilbert spaces to a master equation formulation of stochastic neuronal network models. These theorems are valid for processes taking values in Hilbert spaces and by this are able to incorporate spatial structures of the underlying model.

\end{abstract}

\section{Introduction}\label{sec_intro}

The present study is concerned with the derivation and justification of neural field equations from finite size stochastic particle models, i.e., stochastic models for the behaviour of individual neurons distributed in finitely many populations, in terms of mathematically precise probabilistic limit theorems. We illustrate this approach with the example of the Wilson-Cowan equation
\begin{equation}\label{intro_wc_eq}
\tau\,\dot\nu(t,x)=-\nu(t,x)+f\Bigl(\int_D w(x,y)\nu(t,y)\,\dif y+I(t,x)\Bigr)\,.
\end{equation}
We focus on the following two aspects:
\setlength{\leftmargini}{2.5em}
\begin{enumerate}[align=left]
\setlength{\labelwidth}{2.0em}
  \setlength{\labelsep}{0.5em}
  \setlength{\itemsep}{5pt}
  \setlength{\parsep}{0pt}  
\item[(A)] Often one wants to study deterministic equations such as equation \eqref{intro_wc_eq} in order to obtain results on the `behaviour in the mean' of an intrinsically stochastic system. Thus we first discuss limit theorems of the law of large numbers type for the limit of infinitely many particles. These theorems connect the trajectories of the stochastic particle models to the deterministic solution of mean field equations and hence provide a justification studying equation \eqref{intro_wc_eq} in order to infer on the behaviour of the stochastic system.
\item[(B)] Secondly, we aim to characterise the internal noise structure of the complex discrete stochastic models as in the limit of large numbers of neurons the noise is expected to be close to a simpler stochastic process. Ultimately, this yields a stochastic neural field model in terms of a stochastic evolution equation conceptually analogous to the \emph{Chemical Langevin Equation}. The Chemical Langevin Equation is widely used in the study of chemical reactions networks for which the stochastic effects cannot be neglected but a numerical or analytical study of the exact discrete model is not possible due to its inherent complexity.
\end{enumerate}
In this study we understand as a \emph{microscopic model} a description as a stochastic process, usually a Markov chain model, also called a \emph{master equation formulation} (cf.~\cite{Benayoun,Bressloff1,BuiceCowan,BuiceCowanChow,OhiraCowan2} containing various master equation formulations of neural dynamics). In contrast, a \emph{macroscopic model} is a deterministic evolution equation such as \eqref{intro_wc_eq}. 
Deterministic mean field equations have been used widely and for a long time to model and analyse large scale behaviour of the brain. In their original deterministic form they are successfully used to model geometric visual hallucinations, orientation tuning in the visual cortex and wave propagation in cortical slices to mention only a few applications. We refer to \cite{Bressloff_review} for a recent review and an extensive list of references. 
The derivation of these equations is based on a number of arguments from statistical physics and for a long time a justification from  microscopic models has not been available. 
The interest in deriving mean field equations from stochastic microscopic model has been revived recently as it contains the possibility to derive deterministic `corrections' to the mean field equations, also called second order approximations. 
These corrections might account for the inherent stochasticity and thus incorporate so called finite size effects. This has been achieved by either applying a path-integral approach to the master equation \cite{BuiceCowan,BuiceCowanChow} or by a van Kampen system-size expansion of the master equation \cite{Bressloff1}. 
In more detail, the author in the latter reference proposes a particular master equation for a finite number of neuron populations and derives
the Wilson-Cowan equation as the first order approximation to the mean via employing the van Kampen system size expansion and then taking the continuum limit for a continuum of populations. In keeping also the second order terms a `stochastic' version of the mean field equation is also presented in the sense of coupling the first moment equation to an equation for the second moments.

%
%
%
%

However, the van Kampen system size expansion does not give a precise mathematical connection, as it neither quantifies the type of convergence (quality of the limit), states conditions when the convergence is valid nor does it allow to characterise the speed of convergence.  Furthermore, particular care has to be taken in systems possessing multiple fixed points of the macroscopic equation and we refer to \cite{Bressloff1} for a discussion of this aspect in the neural field setting. The limited applicability of the van Kampen system size expansion was already well known to van Kampen \cite[Sec.~10]{van_Kampen}. 
In parallel to the work of van Kampen, T.~Kurtz derived precise limit theorems connecting sequences of continuous time Markov chains to solutions of systems of ordinary differential equations, see the seminal studies \cite{Kurtz1,Kurtz2} or the monograph \cite{EthierKurtz}. Limit theorems of that type are usually called the \emph{fluid limit, thermodynamic limit} or \emph{hydrodynamic limit}, for a review, see, e.g., \cite{DarlingNorris}.

As is thoroughly discussed in \cite{Bressloff1} establishing the connection between master equation models and mean field equations involves two limit procedures. First, a limit which takes the number of particles, in this case neurons per considered population, to infinity (thermodynamic limit), and a second which gives the mean field by taking the number of populations to infinity (continuum limit). In this `double limit' the theorems by Kurtz describe the connection of taking the number of neurons per population to infinity yielding a system of ordinary differential equation, one for each population. Then the extension from finite to infinite dimensional state space is obtained by a continuum limit. This procedure corresponds to the approach in \cite{Bressloff1}. 
Thus taking the double limit step by step raises the question what happens if we first take the spatial limit and then the fluid limit, thus reversing the order of the limit procedures, or in the case of taking the limits simultaneously. 
Recently, in an extension to the work of Kurtz one of the present authors and co-authors established limit theorems that achieve this double limit \cite{RTW}, thus being able to connect directly finite population master equation formulations to spatio-temporal limit systems, e.g., partial differential equation or integro-differential equations such as the Wilson-Cowan equation \eqref{intro_wc_eq}. In a general framework these limit theorems were derived for Piecewise Deterministic Markov Processes on Hilbert spaces which in addition to the jump evolution also allow for a coupled deterministic continuous evolution. This generality was motivated by applications to neuron membrane models consisting of microscopic models of the ion channels coupled to a deterministic equation for the transmembrane potential. We find that this generality is also advantageous for the present situation of a pure jump model as it allows to include time-dependent inputs. In this study we employ these theorems to achieve the aims (A) and (B) focussing on the example of the deterministic limit given by the Wilson-Cowan equation \eqref{intro_wc_eq}. 

Finally, we state what this study does \emph{not} contain, which in particular distinguishes the present study from \cite{BuiceCowan,BuiceCowanChow,Bressloff1} beyond mathematical technique. Presently, the aim is not to derive moment equations, i.e., a deterministic set of equations that approximate the moments of the Markovian particle model, but rather processes (deterministic or stochastic) to which a sequence of microscopic models converges under suitable conditions in a probabilistic way. 
This means that a microscopic model, which is close to the limit -- presently corresponding to a large number of neurons in a large number of populations --, can be assumed to be close to the limiting processes in structure and pathwise dynamics as indicated by the quality of the stochastic limit. 
Hence, the present work is conceptually -- though neither in technique nor results -- close to \cite{Touboul} wherein using a propagation to chaos approach in the vicinity of neural field equations the author also derives in a mathematically precise way a limiting process to finite particle models. 
However, it is an obvious consequence that the convergence of the models necessarily implies a close resemblance of their moment equations. This provides the connection to \cite{BuiceCowan,BuiceCowanChow,Bressloff1} which we briefly comment on in Appendix \ref{app_moment_equations}.

\medskip

%
%
%
%

As a guide we close this introduction with an outline of the subsequent sections and some general remarks on the notation employed in this study. 
In Sections \ref{section_det_limit} to \ref{sec_inhomogeneou_model_def} we first discuss the two types of mean field models in more detail, on the one hand, the Wilson-Cowan equation as the macroscopic limit and, on the other hand, a master equation formulation of a stochastic neural field.
The main results of the paper are found in Section \ref{sec_precise_LLN}. There we set up the sequence of microscopic models and state conditions for convergence. Limit theorems of the law of large numbers type are presented in Theorem \ref{lln} and Theorem \ref{theorem_infinite_time_lln} in Section \ref{section_llns}. The first is a classical weak law of large numbers providing uniform convergence on compacts in probability and the second convergence in the mean uniformly over the whole positive time axis. Next, a central limit theorem for the martingale part of the microscopic models is presented in Section \ref{section_mart_clt} characterising the internal fluctuations of the model to be of a diffusive nature in the limit. This part of the study is concluded in Section \ref{section_langevin} by presenting the Langevin approximations that arise as a result of the preceding limit theorems.
The proofs of the theorems in Section \ref{sec_precise_LLN} are deferred to Section \ref{section_proofs}. 
The study is concluded in Section \ref{sec_conclusions} with a discussion of the implications of the presented results and an extension of these limit theorems to different master equation formulations or mean field equations.

\medskip 
 
%
%
%
%

\emph{Notations and conventions:} Throughout the study we denote by $L^p(D)$, $1\leq p\leq \infty$, the Lebesgue spaces of real functions on a domain $D\subset\rr^d$, $d\geq 1$. Physically reasonable choices are $d\in\{1,2,3\}$, however for the mathematical theory presented the spatial dimension can be arbitrary. In the present study spatial domains $D$ are always bounded with a sufficiently smooth boundary,  where the minimal assumption is a strong local Lipschitz condition, see \cite{Adams}.  
For bounded domains $D$ this condition simply means that for every point on the boundary its neighbourhood on the boundary is the graph of a Lipschitz continuous function. Furthermore, for $\alpha\in\mathbb{N}$ we denote by $H^\alpha(D)$ the Sobolev spaces, i.e., subspaces of $L^2(D)$, with the corresponding Sobolev norm. For $\alpha\in\rr_+\backslash\mathbb{N}$ we denote by $H^\alpha(D)$ the interpolating Besov spaces. In this study $H^{-\alpha}(D)$ is the dual space of $H^\alpha(D)$ which is in contrast to the widespread notation to denote by $H^{-\alpha}(D)$, $\alpha\geq 0$, the dual space of $H^\alpha_0(D)$. As usual we have $H^0(D)=L^2(D)=H^{-0}(D)$. We thus obtain a continuous scale of Hilbert spaces $H^\alpha(D)$, $\alpha\in\rr$, which satisfy that $H^{\alpha_1}(D)$ is continuously embedded\footnote{A normed space $X$ is continuously embedded in another normed space $Y$, in symbols $X\hookrightarrow Y$, if $X\subset Y$ and there exists a constant $K<\infty$ such that $\|u\|_Y\leq K\|u\|_X$ for all $u\in X$.} in $H^{\alpha_2}(D)$ for all $\alpha_1<\alpha_2$. Next, a pairing $(\,\cdot\,,\,\cdot\,)_{H^\alpha}$ denotes the inner product of the Hilbert space $H^\alpha(D)$ and pairings in angle brackets $\langle\cdot,\cdot\rangle_{H^\alpha}$ denote the duality pairing for the Hilbert space $H^\alpha(D)$. That is, for $\psi\in H^\alpha(D)$ and $\phi\in H^{-\alpha}(D)$ the expression $\langle\phi,\psi\rangle_{H^\alpha}$ denotes the application of the real, linear functional $\phi$ to $\psi$. Furthermore the spaces $H^\alpha(D), L^2(D)$ and $H^{-\alpha}(D)$ form an evolution triplet, i.e., the embeddings are dense and the application of linear functionals and the inner product in $L^2(D)$ satisfy the relation
\begin{equation}\label{evolution_triple_prop}
\langle\phi,\psi\rangle_{H^\alpha}=(\phi,\psi)_{L^2} \quad\forall \phi\in L^2(D), \psi\in H^\alpha(D)\,. 
\end{equation}
Norms in Hilbert spaces are denoted by $\|\cdot\|_{H^\alpha}$, $\|\cdot\|_0$ is used to denote the supremum norm of real functions, i.e., for $f:\rr\to\rr$ we have $\|f\|_0=\sup_{z\in\rr}|f(z)|$, and $|\cdot|$ denotes either the absolute value for scalars or the Lebesgue measure for measurable subsets of Euclidean space. Finally, we use $\mathbb{N}_0$ to denote the set of integers including zero.

\subsection{The macroscopic limit}\label{section_det_limit}

Neural field equations are usually classified into two types, \emph{rate-based} and \emph{activity-based} models. The prototype of the former is the Wilson-Cowan equation, see equation \eqref{intro_wc_eq} which we also restate below, and the Amari equation, see equation \eqref{Amari_equation} in Section \ref{sec_conclusions}, is the prototype of the latter. Besides being of a different structure, due to their derivation, the variable they describe has a completely different interpretation. In rate-based models the variable describes the average rate of activity at a certain location and time, roughly corresponding to the fraction of active neurons at a certain infinitesimal area. In activity-based models the macroscopic variable is an average electrical potential produced by neurons at a certain location. For a concise physical derivation that leads to these models we refer to \cite{Bressloff1}. In the following we consider rate-based equations, in particular, the classical Wilson-Cowan equation, to discuss the type of limit theorems we are able to obtain. We remark, that the results are essentially analogous for activity based models. 

Thus, the macroscopic model of interest is given by the equation
\begin{equation}\label{Wilson_Cowan}
\tau\,\dot\nu(t,x)=-\nu(t,x)+f\Bigl(\int_D w(x,y)\nu(t,y)\,\dif y+I(t,x)\Bigr)\,,
\end{equation}
where $\tau>0$ is a decay time constant, $f:\rr\to\rr_+$ is a gain (or response) function that relates inputs that a neuron receives to activity. In \eqref{Wilson_Cowan} the value $f(z)$ can be interpreted as the fraction of neurons that receive at least threshold input. Furthermore $w(x,y)$ is a weight function which states the connectivity strength of a neuron located at $y$ to a neuron located at $x$ and, finally, $I(t,x)$ is an external input which is received by a neuron at $x$ at time $t$. For the weight function $w:D\times D\to\rr$ and the external input $I$ we assume that $w\in L^2(D\times D)$ and $I\in C(\rr_+,L^2(D))$. As for the gain function $f$ we assume in this study that $f$ is non-negative, satisfies a global Lipschitz condition with constant $L>0$, i.e.,
\begin{equation}\label{lipschitz_on_f}
|f(a)-f(b)|\leq L\,|a-b|\qquad\forall\, a,b\in\rr\,,
\end{equation}
and it is bounded. From an interpretive point-of-view it is reasonable and consistent to stipulate that $f$ is bounded by one -- being a fraction -- as well as being monotone. The latter property corresponds to the fact that higher input results in higher activity. In specific models, $f$ is often chosen to be a sigmoidal function, e.g., $f(z)=(1+\e^{-(\beta_1z+\beta_2})^{-1}$ in \cite{Bressloff2} or $f(z)=(\tanh(\beta_1 z+\beta_2)+1)/2$ in \cite{Benayoun} which both satisfy $f\in[0,1]$. Moreover, the most common choices of $f$ are even infinitely often differentiable with bounded derivatives, which already implies the Lipschitz condition \eqref{lipschitz_on_f}.

The Wilson-Cowan equation \eqref{Wilson_Cowan} is well-posed in the strong sense as an integral equation in $L^2(D)$  under the above conditions. That is, equation \eqref{Wilson_Cowan} possesses a unique, continuously differentiable global solution $\nu$ to every initial condition $\nu(0)=\nu_0\in L^2(D)$, i.e., $\nu\in C^1([0,T],L^2(D))$ for all $T>0$, which depends continuously on the initial condition. Furthermore, if the initial condition satisfies $\nu_0(x)\in [0,\|f\|_0]$ almost everywhere in $D$, then it holds for all $t>0$ that $\nu(t,x)\in (0,\|f\|_0)$ for almost all $x\in D$. For a brief derivation of these results we refer to Section \ref{section_well_posedness} where we also state a result about higher spatial regularity of the solution: Let $\alpha\in\mathbb{N}$ be such that $\alpha>d/2$. If now $\nu_0\in H^\alpha(D)$ and if $f$ is at least $\alpha$-times differentiable with bounded derivatives and the weights and the input function satisfy $w\in H^\alpha(D\times D)$ and $I\in C(\rr_+,H^\alpha(D))$, then the equation is well-posed in $H^\alpha(D)$, i.e., for all $T>0$ in $\nu\in C^1([0,T],H^\alpha(D))$. In particular this implies that the solution $\nu$ is jointly continuous on $\rr_+\times D$.

\subsection{Master equation formulations of neural network models}\label{sec_master_equations}

For the microscopic model we concentrate on a variation of the model considered in \cite{Bressloff1,Bressloff2}, which is already an improvement on a model introduced in \cite{CowanNIPS}. We extend the model including variations among neuron populations and foremost time-dependent inputs. We chose this model over the master equation formulations in \cite{BuiceCowan,BuiceCowanChow} as it provides a more direct connection of the microscopic and macroscopic models, see also the discussion in Section \ref{sec_conclusions}. We describe the main ingredients of the model beginning with the simpler, time-independent model as prevalent in the literature. Subsequently, in Section \ref{sec_inhomogeneou_model_def} the final, time-dependent model is defined.

We denote by $P$ the number of neuron populations in the model. Further, we assume that the $k$-th neuron population consists of identical neurons which can either be in one of two possible states, \emph{active}, i.e., emitting action potentials, and \emph{inactive}, i.e., quiescent or not emitting action potentials. Transitions between states occur instantaneously and at random times. For all $k=1,\ldots,P$ the random variables $\Theta^k_t$ denote the number of active neurons at time $t$. An integer $l(k)$ is used to characterise the population size. This number $l(k)$ can be be interpreted as the number of neurons in the $k$-th population, at least for sufficiently large values. However, this is not accurate in the literal sense as it is possible with positive probability for populations to contain more than $l(k)$ active neurons. Nevertheless, a-posteriori the interpretation can be salvaged from the obtained limit theorems.\footnote{\label{the_bounded_state_space_footnote}The derivation of limit theorems for bounded populations sizes, where $l(k)$ actually \emph{is} the number of neurons per population, is much more delicate than the subsequent presentation as the transition rate functions become discontinuous. Although this would be a desirable result we have not yet been able to prove such a theorem, though it is clear that the Wilson-Cowan equation would be the only possible limit. See also a discussion of this aspect in Section \ref{section_discussion_bounded_state_space}.} It is a corollary of these that the probability of more then $l(k)$ neurons being active for some time becomes arbitrarily small for large enough $l(k)$. Hence, for physiological reasonable neuron numbers the probability in this models of observing `non-physiological' trajectories in the interpretation becomes ever smaller.

Proceeding with notation, $\Theta_t=(\Theta^1_t,\ldots, \Theta^P_t)$ is a (unbounded) piecewise constant stochastic process taking values in $\mathbb{N}_0^P$. The stochastic transitions from inactive to active states and vice versa for a neuron in population $k$ are governed by a constant inactivation rate $\tau^{-1}>0$ -- uniformly for all populations -- and inputs from other neurons depending on the current network state. This non-negative activation rate is given by $\tau^{-1}l(k)\nbar f_k(\theta)$ for $\theta\in \mathbb{N}_0^P$. 
For the definition of $\nbar f_k$ we consider weights $\nbar W_{kj}$, $k,j=1,\ldots, P$, which weight the input one neuron in population $k$ receives from a neuron in population $j$. 
Then the activation rate of a neuron in population $k$ is proportional to
\begin{equation}\label{activation_rate_time_indep}
\nbar f_k(\theta)=f\Bigl(\sum_{j=1}^P \nbar W_{kj}\,\theta^j\Bigr)
\end{equation}
for a non-negative function $f:\rr\to\rr$, which obviously corresponds to the gain function $f$ in the Wilson-Cowan equation \eqref{Wilson_Cowan}. 
We remark that here $f$ is \emph{not} the rate of activation of one neuron. In this model the activation rate of a population is not proportional to the number of inactive neurons but it is proportional to $l(k)$, which stands for the total number of neurons in the population. In \cite{Bressloff1} this rate is thus interpreted as the rate with which a neuron \emph{becomes or remains} active.

It follows that the process $(\Theta_t)_{t\geq 0}$ is a continuous-time Markov chain whose evolution is governed by the following master equation, where $e_k$ denotes the $k$--th basis vector of $\rr^P$,
\begin{eqnarray}
\frac{\dif \Pr[\theta,t]}{\dif t}&=&\frac{1}{\tau}\sum_{k=1}^P\biggl(l(k)\,\nbar f_k(\theta-e_k)\,\Pr[\theta-e_k,t]- \bigl(\theta^k+l(k)\,\nbar f_k(\theta)\bigr)\,\Pr[\theta,t]+ (\theta^k+1)\,\Pr[\theta+e_k,t]\biggr)\phantom{xxx} \label{master_equation_2}
\end{eqnarray}
which is endowed with the boundary conditions $\Pr[\theta,t]=0$ if $\theta\notin \mathbb{N}_0^P$. In \eqref{master_equation_2} the variable $\Pr[\theta,t]$ denotes the probability that the process $\Theta_t$ is in state $\theta$ at time $t$. Finally, the definition is completed with stating an initial law $\mathcal{L}$, the distribution of $\Theta_0$, i.e., providing an initial value for the ODE system \eqref{master_equation_2}.

Another definition of a continuous-time Markov chain is via its generator, see, e.g., \cite{EthierKurtz}, and it is equivalent to the master equation \eqref{master_equation_2}. Although the master equation is widely used in the physics and chemical reactions literature the mathematically more appropriate object for the study of a Markov process is its generator and the master equation is an object derived from the generator, see \cite[Sec.~V]{van_Kampen}. The generator of a Markov process is an operator defined on the space of real functions over the state space of the process. For the above model defined by the master equation \eqref{master_equation_2} the generator is given by
\begin{equation}\label{definition_generator}
\mathcal{A}g(\theta)=\lambda(\theta)\int_{\mathbb{N}^P_0}\Bigl( g(\xi)-g(\theta)\Bigr)\,\mu(\theta,\dif \xi)
\end{equation}
for all suitable $g:\mathbb{N}^P_0\to\rr$. For details we refer to \cite{EthierKurtz}. Here, $\lambda$ is the total instantaneous jump rate, given by
\begin{equation}\label{definition_of_lambda}
\lambda(\theta):= \frac{1}{\tau}\sum_{k=1}^P \Bigl(\theta^k + l(k)\,\nbar f_k(\theta)\Bigr)\,,
\end{equation}
and defines the distribution of the waiting time until the next jump, i.e.,
\begin{equation*}
\Pr[\Theta_{t+s}= \Theta_t\,\forall\, s\in [0,\Delta t]\,|\, \Theta_t=\theta]\,=\, \e^{-\lambda(\theta)\,\Delta t}\,.
\end{equation*}
Further, the measure $\mu$ in \eqref{definition_generator} is a Markov kernel on the state space of the process defining the conditional distribution of the post-jump value, i.e.,
\begin{equation}\label{definition_of_mu}
\Pr[\Theta_t\in A\,|\, \Theta_t\neq \Theta_{t-}]\,=\,\mu(\Theta_{t-},A)
\end{equation}
for all sets $A\subseteq \mathbb{N}_0^P$. In the present case for each $\theta$ the measure $\mu$ is given by the discrete distribution
\begin{equation}\label{definition_of_point_probs}
\mu(\theta,\{\theta-e_k\})=\frac{1}{\tau}\,\frac{\theta^k}{\lambda(\theta)},\quad\mu(\theta,\{\theta+e_k\})=\frac{1}{\tau}\,\frac{l(k)\,\nbar f_k(\theta)}{\lambda(\theta)}\qquad\forall\,k=1,\ldots,P\,.
\end{equation}
The importance of the generator lies in the fact that it fully characterises a Markov process and that convergence of Markov processes is strongly connected to the convergence of their generators, see \cite{EthierKurtz}.

\subsection{Including external time-dependent input}\label{sec_inhomogeneou_model_def}

Until now the microscopic model does not incorporate any time-dependent input into the system. 
In analogy to the macroscopic equation \eqref{Wilson_Cowan} this input enters into the model inside the active rate function $\nbar f_k$. Thus let $\nbar I_k(t)$ denote the external input into a neuron in population $k$ at time $t$, then the time-dependent activation rate is given by
\begin{equation}\label{activation_rate_time_dep}
\nbar f_k(\theta,t)=f\Bigl(\sum_{j=1}^P \nbar W_{kj}\,\theta^j+\nbar I_k(t)\Bigr)\,.
\end{equation}
The most important qualitative difference when substituting \eqref{activation_rate_time_indep} by \eqref{activation_rate_time_dep} is that the corresponding Markov process is no longer homogeneous. In particular the waiting time distributions in between jumps are no longer exponential but satisfy
\begin{equation*}
\Pr[\Theta_{t+s}= \Theta_t\,\forall\, s\in [0,\Delta t]\,|\, \Theta_t=\theta]\,=\, \e^{-\textstyle\int_0^{\Delta t}\lambda(\theta,s)\,\dif s}\,.
\end{equation*}
Hence, the resulting process is an inhomogeneous continuous-time Markov chain, see, e.g., \cite[Sec.~2]{YinZhang}. It is straight forward to write down the corresponding master equation analogously to \eqref{master_equation_2} yielding a system of non-autonomous ordinary differential equations, cf.~the master equation formulation in \cite{BuiceCowan}. Similarly there exists the notion of a time-dependent generator for inhomogeneous Markov processes, cf.~\cite[Sec.~4.7]{EthierKurtz}. 
Employing a standard trick, that is, suitably extending the state space of the process, we can transform a inhomogeneous to a homogeneous Markov process \cite{EthierKurtz,Sharpe}. That is, the space-time process $Y_t:=(\Theta_t,t)$ is again a homogeneous Markov process. The initial law of the associated space-time process is $\mathcal{L}\times\delta_0$ on $\mathbb{N}^P\times\rr_+$. We emphasise that definitions of the space-time process and its initial law imply that the time-component starts at $0$ a.s.~and, moreover, moves continuously and deterministically. That is, the trajectories satisfy in between jumps the differential equation
\begin{equation*}
\left(\begin{array}{c}\dot\theta\\ \dot t\end{array}\right)\,=\, \left(\begin{array}{c} 0 \\ 1\end{array}\right)\,,
\end{equation*}
where the jump intensity $\lambda$ is given by the sum of all individual time-dependent rates analogously to \eqref{definition_of_lambda}. Finally, the post jump value is given by a Markov kernel $\mu((\theta,t),\cdot)\times\delta_t$ as there clearly do not occur jumps in the progression of time and $\mu$ is the obvious time-dependent modification of \eqref{definition_of_point_probs}.

It thus follows, that the space-time process $(\Theta_t,t)_{t\geq 0}$ is a homogeneous Piecewise Deterministic Markov Process (PDMP), see, e.g., \cite{Davis2,Jacobsen,RiedlerPhD}. This connection is particularly important as we apply in the course of the present study limit theorems developed for this type of processes, see \cite{RTW}. Finally, for the space-time process $(\Theta_t,t)_{t\geq 0}$ we obtain for suitable functions $g:\mathbb{N}^P_0\times\rr_+\to\rr$ the generator
\begin{equation}
\mathcal{A}g(\theta,t)=\nabla_{\!t} g(\theta,t)+\lambda(\theta,t)\int_{\mathbb{N}_0^P} \Bigl(g(\xi,t)-g(\theta,t)\Bigr)\,\mu\bigl((\theta,t),\dif \xi\bigr)\,.
\end{equation}

\section{A precise formulation of the limit theorems}\label{sec_precise_LLN}

In this section we present the precise formulations of the limit theorems. To this end we first define a suitable sequence of microscopic models which gives the connection between the defining objects of the Wilson-Cowan equation \eqref{Wilson_Cowan} and the microscopic models discussed in Section \ref{sec_master_equations}. Thus, $(Y_t^n)_{t\geq 0}=(\Theta^n_t,t)_{t\geq 0}$, $n\in\mathbb{N}$, denotes a sequence of microscopic PDMP neural field models of the type as defined in Section \ref{sec_inhomogeneou_model_def}. Each process $(Y_t^n)_{t\geq 0}$ is defined on a filtered probability space $(\Omega^n,\sF^n,(\sF^n_t)_{t\geq 0},\Pr^n)$ which satisfies the usual conditions. Hence, the defining objects for the jump models are now dependent on an additional index $n$. That is $P(n)$ denotes the number of neuron populations in the \mbox{$n$-th} model, $l(k,n)$ is the number of neurons in the $k$-th population of the \mbox{$n$-th} model and analogously we use the notations $\nbar W_{kj}^n$ and $\nbar I_{k,n}$ and $\nbar f_{k,n}$. However, we note from the beginning that the decay rate $\tau^{-1}$ is independent of $n$ and $\tau$ is the time constant in the Wilson-Cowan equation \eqref{Wilson_Cowan}. In the following paragraphs we discuss the connection of the defining components of this sequence of microscopic models to the components of the macroscopic limit.\medskip

\textbf{Connection to the spatial domain $D$.} A key step of connecting the microscopic models to the solution of equation \eqref{Wilson_Cowan} is that we need to put the individual neuron populations into relation to the spatial domain $D$ the solution of \eqref{Wilson_Cowan} lives on. To this end we assume that each population is located within a subdomain of $D$ and that the subdomains of the individual populations are non-overlapping. Hence, for each $n\in\mathbb{N}$ we obtain a collection $\mathcal{D}_n$ of $P(n)$ non-overlapping subsets of $D$ denoted by $D_{1,n},\ldots, D_{P(n),n}$. We assume that each subdomain is measurable and convex. The convexity of the subdomains is a technical condition that allows us to apply Poincar\'e's inequality, cf.~\eqref{poincare_inequality}. We do not think that this condition is too restrictive as most reasonable partition domains, e.g., cubes, triangles, are convex. Furthermore, for all reasonable domains $D$, e.g., all Jordan measurable domains, a sequence of convex partitions can be found such that additionally the conditions imposed in the limit theorems below are also satisfied. Conversely, one may think of obtaining the collection  $\mathcal{D}_n$ by partitioning the domain into $P(n)$ convex subdomains $D_{1,n},\ldots, D_{P(n),n}$ and confining each neuron population to one subdomain. However it is not required that the union of the sets in $\mathcal{D}_n$ amounts to the full domain $D$ nor that the partitions consists of refinements. Necessary conditions on the limiting behaviour of the subdomains are very strongly connected to the convergence of initial conditions of the models, which is a condition in the limit theorems, see below. For the sake of terminological simplicity we refer to $\mathcal{D}_n$ simply as the partitions.

We now define some notation for parameters characterising the partitions $\mathcal{D}_n$: the minimum and maximum Lebesgue measure, i.e., length, area or volume depending on the spatial dimension, is denoted by
\begin{equation}
v_-(n):=\min_{k=1,\ldots,P(n)} |D_{k,n}|,\qquad v_+(n):=\max_{k=1,\ldots,P(n)} |D_{k,n}|\,,
\end{equation}
and the maximum diameter of the partition is denoted by
\begin{equation}
\delta_+(n):=\max_{1,\ldots,P(n)}\,\textnormal{diam}\,(D_{k,n}) \,,
\end{equation}
where the diameter of a set $D_{k,n}$ is defined as $\textnormal{diam}\,(D_{k,n}):=\sup_{x,y\in D_{k,n}}|x-y|$. In the special case of domains obtained by unions of cubes with edge length $n^{-1}$ it obviously holds that $v_\pm(n)=n^{-d}$ and $\delta_+(n)=\sqrt{d}\,n^{-1}$. It is a necessary condition in all the limit theorems that $\lim_{n\to\infty}\delta_+(n) = 0$ which implies that $\lim_{n\to\infty}v_+(n) = 0$ as well as $\lim_{n\to\infty} P(n) = \infty$ as the Lebesgue measure of a set is bounded in terms of the diameter of the set. That is, in order to obtain a limit the sequence of partitions necessarily consists of ever finer sets and the number of neuron populations has to diverge. Finally, each domain $D_{k,n}$ of the partition $\mathcal{D}_n$ contains one neuron population `consisting' of $l(k,n)\in\mathbb{N}$ neurons. Then we denote by $\ell_\pm(n)$ the maximum and minimum number of neurons in populations corresponding to the $n$-th model, i.e.,
\begin{equation}
\ell_-(n):=\min_{k=1,\ldots,P(n)} l(k,n),\qquad \ell_+(n):=\max_{k=1,\ldots,P(n)} l(k,n)\,.
\end{equation}

\medskip

\textbf{Connection to the weight function $w$.} We assume that there exists a function $w: D\times D \to \rr$ such that the connection to the discrete weights is given by
\begin{equation}\label{definition_of_discret_weights}
\nbar W_{kj}^{n}:=\frac{1}{|D_{k,n}|}\int_{D_{k,n}}\biggl(\int_{D_{j,n}}w(x,y)\,\dif y\biggr)\,\dif x\,,
\end{equation}
where $w$ is the same function as in the Wilson-Cowan equation \eqref{Wilson_Cowan}. For the definition of activation rate at time $t$ we thus obtain
\begin{equation}\label{definition_of_discret_increase_rates}
\nbar f_{k,n}(\theta^n\!,t):=f\biggl(\sum_{j=1}^{P} \nbar W^{n}_{kj}\frac{\theta^{j,n}}{l(j,n)}+\nbar I_{k,n}(t)\biggr)\,. 
\end{equation}

\medskip

\textbf{Connection to the input current $I$.} The external input which is applied to neurons in a certain population is obtained by spatially averaging a space-time input over the subdomain that population is located in, i.e.,
\begin{equation}\label{definition_of_discret_inputs}
\nbar I_{k,n}(t):= \frac{1}{|D_{k,n}|}\int_{D_{k,n}} I(t,x)\dif x\,.
\end{equation}

This completes the definition of the Markov jump processes $(\Theta^n_t,t)_{t\geq 0}$. For the sake of completeness we repeat the definition of the total jump rate
\begin{equation*}
 \lambda^n(\theta^n\!,t):=\frac{1}{\tau}\sum_{k=1}^P\Bigl(\theta^{k,n}+l(k,n)\,\nbar f_{k,n}(\theta^n\!,t)\Bigr)
\end{equation*}
and the transition measure $\mu^n$ is defined by
\begin{equation*}
\mu^n\bigl((\theta^n\!,t),\{\theta^n-e_k\}\bigr):=\frac{1}{\tau}\frac{\theta^{k,n}}{\lambda^n(\theta^n\!,t)},\qquad \mu^n\bigl((\theta^n\!,t),\{\theta^n+e_k\}\bigr):=\frac{1}{\tau}\frac{l(k,n)\nbar f_{k,n}(\theta^n\!,t)}{\lambda^n(\theta^n\!,t)}
\end{equation*}
for all $k=1,\ldots P(n)$.

\medskip

\textbf{Connection to the solution $\nu$.} As functions of time, the paths of the PDMP $(\Theta^n_t,t)_{t\geq 0}$ and the solution $\nu$ live on different state spaces. The former takes values in $\mathbb{N}^P_0\times \rr_+$ and the latter in $L^2(D)$. Thus in order to compare these two we have to introduce a mapping that maps the stochastic process onto $L^2(D)$. In \cite{RTW} the authors called such a mapping a \emph{coordinate function} which is also the terminology used in \cite{DarlingNorris}. In fact, the limit theorems we subsequently present actually are for the processes we obtain from the composition of the coordinate functions with the PDMPs. Here it is important to note that for each $n\in\mathbb{N}$ the coordinate functions may -- and usually do -- differ, however they project the process into the common space $L^2(D)$. For the mean field models we define the coordinate functions for all $n\in\mathbb{N}$ by 
\begin{equation}\label{def_coordinate_function}
\nu^n:\mathbb{N}_0^P\to L^2(D): \theta^n\mapsto\sum_{k=1}^P\frac{\theta^{k,n}}{l(k,n)}\,\mathbb{I}_{D_{k,n}}\,.
\end{equation}
Clearly each $\nu^n$ is a measurable map into $L^2(D)$. For the composition of $\nu^n$ with the stochastic process $(\Theta^n_t,t)_{t\geq 0}$ we also use the abbreviation $\nu^n_t:=\nu^n(\Theta^n_t)$ and hence the resulting stochastic process $(\nu^n_t)_{t\geq 0}$ is an adapted c\`adl\`ag process taking values in $L^2(D)$. This process thus states the activity at a location $x\in D$ as the fraction of active neurons in the population which is located around this location.

\medskip

\textbf{Connection of the initial conditions.} One condition in the subsequent limit theorems is the convergence of initial conditions in probability, i.e., the assumption that
\begin{equation}\label{prob_conv_initial_cond}
\lim_{n\to\infty} \Pr^n\bigl[\|\nu^n(\Theta_0^n)-\nu_0\|_{L^2}>\eps\bigr]=0\qquad\forall\, \eps>0\,.
\end{equation}
It is easy to see that such a sequence of initial conditions $\Theta^n_0$, $n\in\mathbb{N}$, can be found for any deterministic initial condition $\nu_0$ under some reasonable conditions on the domain $D$ and the sequence of partitions $\mathcal{D}_n$. Hence the assumption \eqref{prob_conv_initial_cond} can always be satisfied. For example, we may define such a sequence of initial conditions by
\begin{equation*}
\Theta^{k,n}_0=\textnormal{argmin}_{i=1,\ldots,l(k,n)}\Big|\frac{i}{l(k,n)}-\frac{1}{|D_{k,n}|}\int_{D_{k,n}}\nu(0,x)\,\dif x\Big|\,.
\end{equation*}
Next, assuming that partitions fill the whole domain $D$ for $n\to\infty$, i.e., $\lim_{n\to\infty} |D\backslash\bigcup_{k=1}^{P(n)}D_{k,n}|=0$, and that the maximal diameter of the sets decreases to zero, i.e., $\lim_{n\to\infty}\delta_+(n)=0$, it is easy to see using the Poincar\'e inequality \eqref{poincare_inequality} that the above definition of the initial condition implies that $\|\nu^n_0-\nu(0)\|_{L_2}\to 0$ and $\sup_{n\in\mathbb{N}}\|\nu_0^n\|_{L^2}^{2r}<\infty$ for all $r\geq 1$. Then \eqref{prob_conv_initial_cond} holds trivially as the initial condition is deterministic and converges. 
A simple non-degenerate sequence of initial conditions is obtained by choosing random initial conditions with the above value as their mean and sufficiently fast decreasing fluctuations. 
Furthermore, a sequence of partitions which satisfy the above conditions also exists for a large class of reasonable domains $D$. Assume that $D$ is Jordan measurable, i.e., a bounded domain such that the boundary is a Lebesgue null set, and let $\mathcal{C}_n$ be the smallest grid of cubes with edge length $1/n$ covering $D$. We define $\mathcal{D}_n$ to be the set of all cubes which are fully in $D$. As $D$ is Jordan measurable these partitions fill up $D$ from inside and $\delta_+(n)\to 0$. For a more detailed discussion of these aspects we refer to \cite{RiedlerPhD}.

\medskip

In the remainder of this section we now collect the main results of this article. We start with the law of large numbers, which establishes the connection to the deterministic mean field equation, and then proceed to central limit theorems which provide the basis for a Langevin approximation. The proofs of the results are deferred to Section \ref{section_proofs}.

\subsection{A law of large numbers}\label{section_llns}

The first law of large numbers takes the following form. Note that the assumptions imply that the number of neuron populations diverges.

\begin{theorem}[Law of large numbers]\label{lln} Let $w\in L^2(D)\times L^2(D)$ and $I\in L^2_{\textnormal{loc}}(\rr_+,H^1(D))$. Assume that  the sequence of initial conditions converges to $\nu(0)$ in probability in the space $L^2(D)$, i.e., \eqref{prob_conv_initial_cond} holds, that $\EX^n\Theta^{k,n}_0\leq l(k,n)$ and that
\begin{equation}\label{assumptions_lln}
\lim_{n\to\infty}\delta_+(n) = 0,\quad \lim_{n\to\infty} \ell_-(n)=\infty
\end{equation}
holds. Then it follows that the sequence of $L^2(D)$--valued jump-process $(\nu^n_t)_{t\geq 0}$ converges in probability uniformly on compact time intervals to the solution $\nu$ of the Wilson-Cowan equation \eqref{Wilson_Cowan}, i.e., for all $T,\eps>0$ it holds that
\begin{equation}\label{thm_lln_conv_in_prob}
\lim_{n\to\infty}\Pr^n\bigl[\sup\nolimits_{t\in[0,T]}\|\nu^n_t-\nu(t)\|_{L^2(D)}>\eps\bigr]=0\,.
\end{equation}
Moreover, if for $r\geq 1$ the initial conditions satisfy in addition $\sup_{n\in\mathbb{N}}\EX^n\|\nu^n_0\|_{L^2}^{2r}<\infty$, then convergence in the $r$-th mean holds, i.e., for all $T>0$
\begin{equation}\label{thm_lln_conv_in_mean}
\lim_{n\to\infty} \EX^n\sup\nolimits_{t\in[0,T]}\|\nu^n_t-\nu(t)\|^r_{L^2(D)}=0\,.
\end{equation}
\end{theorem}

\begin{remark}
The norm of the uniform convergence $\sup_{t\in[0,T]}\|\cdot\|_{L^2}$ for which we have stated convergence in probability and in the mean in Theorem \ref{lln} is a very strong norm on the space of $L^2(D)$--valued c\`adl\`ag functions on $[0,T]$. Hence, due to continuous embeddings the result immediately extends to weaker norms, e.g., the norms $L^p((0,T),L^2)$ for all $1\leq p\leq \infty$. Also for the state space weaker spatial norms can be chosen, e.g., $L^p(D)$ with $1\leq p\leq 2$ or any norm on the duals $H^{-\alpha}(D)$ of Sobolev spaces with $\alpha>0$. 
If weaker norms for the state space are considered it is even possible to relax the conditions of Theorem \ref{lln} by sharpening some estimates in the proof of the theorem. Clearly, it is sufficient that the initial conditions converge in probability with respect to the weaker norms. Recall that $H^{-\alpha}(D)$ denotes the dual of the Sobolev space $H^\alpha(D)$ and $H^0(D)=L^2(D)=H^{-0}(D)$. The results in the following corollary cover the whole range of $\alpha\geq 0$ and splits it into sections with weakening conditions. In particular note, that after passing to weaker norms the convergence does not necessitate that the neuron numbers per population diverge. However, regarding the divergence of the neuron populations, this condition ($\delta_+(n)\to 0$) cannot be relaxed. 

\end{remark}

\begin{corollary}\label{corollary_to_lln} Let $\alpha\geq 0$ and set
\begin{equation}\label{corollary-def_of_q}
q:=\left\{\begin{array}{cl}
 \frac{2d}{d+2\alpha} &\textnormal{ if } 0\leq \alpha<d/2,\\[1ex] 1- &\textnormal{ if } \alpha=d/2,\\[1ex] 1 &\textnormal{ if } d/2< \alpha <\infty.
\end{array}\right. 
\end{equation}
Further, assume that $w\in L^q(D)\times L^2(D)$ and $I\in L^2_\textnormal{loc}(\rr_+,H^1(D))$ and that the sequence of initial conditions converges to $\nu(0)$ in probability in the space $H^{-\alpha}(D)$, that $\lim_{n\to\infty} \delta_+(n)=0$ and
\begin{equation}\label{lln_sharpened_conditions}
\left.\begin{array}{ll}
\displaystyle\lim_{n\to\infty} \frac{v_+(n)^{2\alpha/d}}{\ell_-(n)}=0& \textnormal{ if } 0\leq \alpha<d/2,\\[2ex] 
\displaystyle\lim_{n\to\infty} \frac{v_+(n)^{1-}}{\ell_-(n)}=0 & \textnormal{ if } \alpha=d/2,\\[2ex]
\displaystyle\lim_{n\to\infty} \frac{v_+(n)}{\ell_-(n)}=0 & \textnormal{ if } d/2< \alpha <\infty\,,
\end{array}\right\}
\end{equation}
where $1-$ denotes an arbitrary positive number strictly smaller than $1$. Then it holds for all $T,\eps>0$ that
\begin{equation*}
\lim_{n\to\infty}\Pr^n\bigl[\sup\nolimits_{t\in[0,T]}\|\nu^n_t-\nu(t)\|_{H^{-\alpha}(D)}>\eps\bigr]=0
\end{equation*}
and for $r\geq 1$, if the additional boundedness assumptions of \emph{Theorem \ref{lln}} are satisfied, that for all $T>0$
\begin{equation*}
\lim_{n\to\infty} \EX^n\sup\nolimits_{t\in[0,T]}\|\nu^n_t-\nu(t)\|^r_{H^{-\alpha}(D)}=0\,.
\end{equation*}
\end{corollary}

\begin{remark} We believe that fruitful and illustrative comparisons of these convergence results and their conditions to the results in Kotelenez \cite{Kotelenez1,Kotelenez2} and, particularly, Blount \cite{Blount1} can be made. Here we just mention that the latter author conjectured the conditions \eqref{lln_sharpened_conditions} to be optimal for the convergence but was not able to prove this result in his model of chemical reactions with diffusions for the region $\alpha\in (0,d/2]$. For our model we could achieve these rates.
\end{remark}

\subsubsection{Infinite-time convergence}

In the law of large numbers, Theorem \ref{lln}, and its Corollary \ref{corollary_to_lln} we have presented results of convergence over finite time intervals. Employing a different technique, we are also able to derive a convergence result over the whole positive time axis motivated by a similar result in \cite{TouboulHermann}. The proof of the following theorem is deferred to Section \ref{section_proof_infinite_time_lln}. Restricted to finite time intervals the subsequent result is strictly weaker than Theorem \ref{lln}. However, the result is important when one wants to analyse the mean long time behaviour of the stochastic model via a bifurcation analysis of the \mbox{deterministic} limit as \eqref{infinite_time_lln_result} suggests that $\EX^n\nu^n_t$ is close to $\nu(t)$ for all times $t\geq 0$ for sufficiently large $n$.

\begin{theorem}\label{theorem_infinite_time_lln} Let $\alpha\geq 0$ and assume that the conditions of \textnormal{Corollary \ref{corollary_to_lln}} are satisfied. We further assume that the current input function $I\in L^2_{\textnormal{loc}}(\rr_+,H^1(D))$ satisfies $\|\nabla_{\!x} I\|_{L^\infty(\rr_+,L^2(D))}<\infty$, i.e., it is square integrable in $H^1(D)$ over bounded intervals, and possesses first spatial derivatives bounded for almost all $t\geq 0$ in $L^2(D)$. Then it holds that
\begin{equation}\label{infinite_time_lln_result}
\lim_{n\to\infty} \sup\nolimits_{t\geq 0}\,\EX^n\|\nu^n_t-\nu(t)\|_{H^{-\alpha}(D)} \,=\,0\,.
\end{equation}
\end{theorem}

\subsection{A martingale central limit theorem}\label{section_mart_clt}

In this section we present a central limit theorem for a sequence of martingales associated with the jump processes $\nu^n$. A brief, heuristic discussion of the method of proof for the law of large numbers explains the importance of these martingales and motivates their study.
In the proof of the law of large numbers the central argument relies on the fact that the process $(\nu^n_t)_{t\geq 0}$ satisfies the decomposition
\begin{equation}\label{expansion_of_jump_process}
\nu^n_t= \nu^n_0 + \int_0^t \lambda(\Theta^n_s,s)\int_{\mathbb{N}^P_0} \bigl(\nu^n(\xi)-\nu^n(\Theta^n_s)\bigr)\,\mu^n\bigl((\Theta^n_s,s),\dif\xi\bigr)\,\dif s + M^n_t\,.
\end{equation}
Here the process $(M^n_t)_{t\geq 0}$ is a Hilbert space-valued, square-integrable, c\`adl\`ag martingale and, using \eqref{expansion_of_jump_process} as its definition, is given by
\begin{equation*}
M^n_t=\nu^n_t-\nu^n_0-\int_0^t \lambda(\Theta^n_s,s)\int_{\mathbb{N}^P_0} \bigl(\nu^n(\xi)-\nu^n(\Theta^n_s)\bigr)\,\mu^n\bigl((\Theta^n_s,s),\dif\xi\bigr)\,\dif s\,.
\end{equation*}
We have also used this representation of the process $\nu^n$ in the proof of Theorem \ref{theorem_infinite_time_lln}, see Section \ref{section_proof_infinite_time_lln}. We note that the Bochner integral in \eqref{expansion_of_jump_process} is a.s.~well defined due to bounded second moments of the integrand, see \eqref{proof_lln_cond_a} in the proof of Theorem \ref{lln}. Now an heuristic argument to obtain the convergence to the solution of the Wilson-Cowan equation is the following: The initial conditions converge, the martingale term $M^n$ converges to zero and the integral term in the right hand side of \eqref{expansion_of_jump_process} converges to the right hand side in the Wilson-Cowan equation \eqref{Wilson_Cowan}. Hence, the `solution' $\nu^n$ of \eqref{expansion_of_jump_process} converges to the solution $\nu$ of the Wilson-Cowan equation \eqref{Wilson_Cowan}. 
Now interpreting equation \eqref{expansion_of_jump_process} as a stochastic evolution equation which is driven by the martingale $(M^n_t)_{t\geq 0}$ sheds light on the importance of the study of this term. Because, from this point of view the martingale part in the decomposition \eqref{expansion_of_jump_process} contains all the stochasticity inherent in the system.
Then the idea for deriving a Langevin or linear noise approximation is to find a stochastic non-trivial limit (in distribution) for the sequence of martingales and substituting heuristically this limiting martingale into the stochastic evolution equation. Then it is expected that this new and much less complex process behaves similarly to the process $(\nu^n_t)_{t\geq 0}$ for sufficiently large $n$. 
Deriving a suitable limit for $(M^n_t)_{t\geq 0}$ is what we set to do next. The result can be found in Theorem \ref{martingale_clt} below and takes the form of a central limit theorem.\medskip

First of all, what has been said so far implies the necessity of rescaling the martingale with a diverging sequence in order to obtain a non-trivial limit. The conditions in the law of large numbers imply in particular that the martingale converges uniformly in the mean square to zero, i.e.,
\begin{equation*}
\lim_{n\to\infty} \EX^n\sup\nolimits_{\,t\in[0,T]}\|M^n_t\|_{L^2}\ = \ 0,
\end{equation*}
which in turn implies convergence in probability and convergence in distribution to the zero limit. 
%


Furthermore, in contrast to Euclidean spaces norms on infinite-dimensional spaces are usually not equivalent. 
In Corollary \ref{corollary_to_lln} we exploited this fact as it allowed us to obtain convergence results under less restrictive conditions by changing to strictly weaker norms. In the formulation and proof of central limit theorems, the change to weaker norms even becomes an essential ingredient. 
It is often observed in the literature, see, e.g., \cite{Blount1,Kotelenez1,Kotelenez2}, that central limit theorems cannot be proven in the strongest norm for which the law of large numbers holds, e.g., $L^2(D)$ in the present setting, but only in a strictly weaker norm. Here this norm is the norm in the dual of an appropriate Sobolev space. Hence, from now on we consider for all $n\in\mathbb{N}$ the processes $(\nu^n_t)_{t\geq 0}$ and the martingales $(M^n_t)_{t\geq 0}$ as taking values in the space $H^{-\alpha}(D)$ for an $\alpha> d$, where $d$ is the dimension of the spatial domain $D$, using the embedding of $L^2(D)$ into $H^{-\alpha}(D)$. The technical significance of the restriction $\alpha>d$ is that these are the indices such that there exists an embedding $H^{\alpha}(D)$ into a $H^{\alpha_1}(D)$ with $d/2<\alpha_1<\alpha$ which is of Hilbert-Schmidt type\footnote{\label{maurin_footnote}A continuous embedding of two Hilbert spaces $X\hookrightarrow Y$ is of Hilbert-Schmidt type if for every orthonormal basis $\varphi_j$, $j\in\mathbb{N}$, of $X$ it holds that $\sum_{j=1}^\infty \|\varphi_j\|_{Y}^2<\infty$. Then, more precisely, Maurin's Theorem states that for non-negative integers $m,k$ the embedding of $H^{m+k}(D)$ into $H^m(D)$ is of Hilbert-Schmidt type for $k>d/2$, see \cite{Adams}. The result was generalized to fractional order Sobolev spaces in \cite{Wloka}: Let $D$ be a bounded, strong local Lipschitz domain in $\rr^d$ and $0\leq \alpha_1<\alpha_2$ are real numbers. Then it holds that the embedding of $H^{\alpha_2+d/2}(D)$ into $H^{\alpha_1}(D)$ is of Hilbert-Schmidt type.} due to Maurin's Theorem and $H^{\alpha_1}(D)$ is embedded into $C(\nbar D)$ due to the Sobolev Embedding Theorem. These two properties are essential for the proof of the central limit theorem and their occurrence will be made clear subsequently. \medskip

The limit we propose for the rescaled martingale sequence is a \emph{centred diffusion process} in $H^{-\alpha}(D)$. That is, a centred Gaussian stochastic process $(X_t)_{t\geq 0}$ taking values in $H^{-\alpha}(D)$ with independent increments and given covariance $C(t)$, $t\geq 0$, see, e.g., \cite{DaPratoZabczyk,PrevotRoeckner} for a discussion of Gaussian processes in Hilbert spaces. Such a process is uniquely defined by its covariance operator and conversely, each family of linear, bounded operators $C(t):H^{\alpha}(D)\to H^{-\alpha}(D)$, $t\geq 0$, uniquely defines a diffusion process\footnote{Usually the covariance operator for a Hilbert space-valued process is an operator mapping from the state space into the state space and not into the dual, i.e., in the present situation mapping $H^{-\alpha}(D)$ into itself. Due to the canonical embedding of Hilbert spaces into their dual and the Riesz Representation, however, we can effortless change from the usual definition to ours and vice versa. Moreover, the symmetry condition thus implies due to the Hellinger-Toeplitz Theorem that the operator is self-adjoint and hence of trace class if and only if \eqref{trace_class_cond} is satisfied. The choice of the presentation here is due to the fact that it is simpler to evaluate the duality pairing on $H^{-\alpha}(D)$ than the inner product thereon, as the former usually is just the inner product in $L^2(D)$.} if
\setlength{\leftmargini}{2.0em}
\begin{enumerate}[align=left]
\setlength{\labelwidth}{1.5em}
  \setlength{\labelsep}{0.5em}
  \setlength{\itemsep}{0pt}
  \setlength{\parsep}{0pt}  
\item[(i)] each $C(t)$ is \emph{symmetric} and \emph{positive}, i.e.,
\begin{equation*}
\langle C(t)\,\phi,\psi\rangle_{H^\alpha}=\langle C(t)\,\psi,\phi\rangle_{H^\alpha}\quad\textnormal{and}\quad \langle C(t)\,\phi,\phi\rangle_{H^\alpha}\geq 0\,,
\end{equation*}
\item[(ii)] each $C(t)$ is of \emph{trace class}, i.e., for one (and thus every) orthonormal basis $\varphi_j$, $j\in\mathbb{N}$, in $H^\alpha(D)$ it holds that
\begin{equation}\label{trace_class_cond}
 \sum_{j=1}^\infty\langle C(t)\varphi_j,\varphi_j\rangle_{H^\alpha}<\infty,
\end{equation}
\item[(iii)] and the family $C(t)$, $t\geq 0$, is \emph{continuously increasing} in $t$ in the sense that the map $t\mapsto \langle C(t)\,\phi,\psi\rangle_{H^\alpha}$ is continuous and increasing for all $\phi,\psi\in H^\alpha(D)$.
\end{enumerate}

We next define the process which will be the limit identified in the martingale central limit theorem via its covariance. In order to define the operator $C$ we first define a family of linear operators $G(\nu(t),t)$ mapping from $H^\alpha(D)$ into the dual space $H^{-\alpha}(D)$ via the bilinear form
\begin{equation}\label{def_quad_var_diffusion}
\langle G(\nu(t),t)\,\phi,\psi\rangle_{H^\alpha} \,=\, \int_D \phi(x)\Bigl(\frac{1}{\tau}\nu(t,x)+\frac{1}{\tau}\,f\Bigl(\int_Dw(x,y)\,\nu(t,y)\,\dif y+I(t,x)\Bigr)\Bigr)\,\psi(x)\,\dif x\,.
\end{equation}
It is obvious that this bilinear form is symmetric and positive and, as $\nu(t)$ is continuous in $t$, it holds that the map $t\mapsto \langle G(\nu(t),t)\,\phi,\psi\rangle_{H^\alpha}$ is continuous for all $\phi,\psi\in H^\alpha(D)$. Furthermore, it is easy to see that the operator is bounded, i.e.,
\begin{equation*}
\|G(\nu(t),t)\|_{L(H^\alpha,H^{-\alpha})}=\sup_{\|\phi\|_{H^\alpha}=1}\|G(\nu(t),t)\,\phi\|_{H^{-\alpha}} = \sup_{\|\phi\|_{H^\alpha}=1}\sup_{\|\psi\|_{H^\alpha}=1}\big|\langle G(\nu(t),t)\,\phi,\psi\rangle_{H^\alpha}\big|  < \infty,
\end{equation*}
as the solution of the Wilson-Cowan equation $\nu$ and the gain function $f$ are pointwise bounded. Hence due to the Cauchy-Schwarz inequality the norm $\big|\langle G(\nu(t),t)\,\phi,\psi\rangle_{H^\alpha}\big|$ is proportional to the product $\|\phi\|_{L^2}\|\psi\|_{L^2}$ and for any $\alpha\geq 0$ the Sobolev Embedding Theorem gives now a uniform bound in terms of the norm of $\phi,\,\psi$ in $H^\alpha(D)$. As a final property we show that these operators are of trace-class if $\alpha>d/2$. Thus let $(\varphi_j)_{j\in\mathbb{N}}$ be an orthonormal basis in $H^{\alpha}(D)$, then the Cauchy-Schwarz inequality yields
\begin{equation*}
\big|\langle G(\nu(t),t)\,\varphi_j,\varphi_j\rangle_{H^\alpha}\big|\,\leq\,\frac{1}{\tau}\,(1+\|f\|_0)\,|D|\,\|\varphi_j\|^2_{L^2}\,.
\end{equation*}
Summing these inequalities for all $j\in\mathbb{N}$ we find that the resulting right hand side is finite as due to Maurin's Theorem the embedding of $H^{\alpha}(D)$ into $L^2(D)$ is of Hilbert-Schmidt type. Moreover their trace is even bounded independently of $t$. 

Now, it holds that the map $t\mapsto G(\nu(t),t)$ is continuous taking values in the Banach space of trace class operators, hence we define trace class operators $C(t)$ from $H^{\alpha}(D)$ into $H^{-\alpha}(D)$ via the Bochner integral for all $t\geq 0$
\begin{equation}\label{def_covariance_op}
C(t):=\int_0^t G(\nu(s),s)\,\dif s\,.
\end{equation} 
Clearly, the resulting bilinear form $\langle C(t)\cdot,\cdot\rangle_{H^{\alpha}}$ inherits the properties of the bilinear form \eqref{def_quad_var_diffusion}. Moreover, due to the positivity of the integrands it follows that $\langle C(t)\phi,\phi\rangle_{H^\alpha}$ is increasing in $t$ for all $\phi\in H^\alpha(D)$. Hence the family of operators $C(t),\,t\geq 0$, satisfies the above conditions (i)--(iii) and thus uniquely defines an $H^{-\alpha}(D)$--valued diffusion process.\medskip

We are now able to state the martingale central limit theorem. The proof of the theorem is deferred to Section \ref{section_proof_martingale_clt}.

\begin{theorem}[Martingale central limit theorem]\label{martingale_clt} Let $\alpha>d$ and assume that the conditions of \textnormal{Theorem \ref{lln}} are satisfied. In particular convergence in the mean holds, i.e., \eqref{thm_lln_conv_in_mean} holds for $r=1$. Additionally, we assume it holds that
\begin{equation}\label{assumption_on_volumes}
\lim_{n\to\infty} \frac{v_-(n)}{v_+(n)}\frac{\ell_-(n)}{\ell_+(n)}=1\,.
\end{equation}
Then it follows that the sequence of rescaled $H^{-\alpha}(D)$--valued martingales
\begin{equation*}
\Bigl(\sqrt{\tfrac{\ell_-(n)}{v_+(n)}}\,M^n_t\Bigr)_{t\geq 0} 
\end{equation*}
converges weakly to the $H^{-\alpha}(D)$--valued diffusion process defined by the covariance operator $C(t)$ given by \eqref{def_covariance_op}.
\end{theorem}

\begin{remark} In connection with the results of Theorem \ref{martingale_clt} two questions may arise. First, in what sense is there uniqueness of the rescaling sequence and hence of the limiting diffusion? That is, does a different scaling also produce a (non-trivial) limit, or, rephrased, is the proposed scaling the correct one to look at? Secondly, the theorem deals with the norms for the range of $\alpha>d$ in the Hilbert scale, what can be said about convergence in the stronger norms corresponding to the range of $\alpha\in[0,d]$? Does there exist a limit? We conclude this section addressing these two issues. 

Regarding the first question, it is immediately obvious that the rescaling sequence $\tfrac{\ell_-(n)}{v_+(n)}$, which we denote by $\rho_n$ in the following, is not a unique sequence yielding a non-trivial limit. Rescaling the martingales $M^n$ by any sequence of the form $ \sqrt{c\rho_n}$ yields a convergent martingale sequence. However, the limiting diffusion differs only in a covariance operator which is also rescaled by $c$ and hence the limit is essentially the same process with either `stretched' or `shrinked' variability. However, the asymptotic behaviour of the rescaling sequences which allow for a non-trivial weak limit is unique.
In general, by considering different rescaling sequences $\rho_n^\ast$ we obtain three possibilities for the convergence of the sequence $\sqrt{\rho^\ast_n}\,M^n$. If $\rho^\ast_n$ is of the same speed of convergence as $\rho_n$, i.e., for $\rho^\ast_n=\landau(\rho_n)$, the thus rescaled sequence converges again to a diffusion process for which the covariance operator is proportional to \eqref{def_covariance_op}. This is then just a rescaling by a sequence (asymptotically) proportional to $\rho_n$ as discussed above. 
Secondly, if the convergence is slower, i.e., $\rho_n^\ast=o(\rho_n)$, then the same methods as in the law of large numbers show that the sequence converges uniformly on compacts in probability to zero, hence also convergence in distribution to the degenerate zero process follows. Thus one only obtains the trivial limit.
Finally, if we rescale by a sequence that diverges faster, i.e., $\rho_n=o(\rho^\ast_n)$, we can show that there does not exist a limit. This follows from general necessary conditions for the preservation of weak limits under transformation which presuppose that $\sqrt{\rho_n^\ast/\rho_n}\,M$ has to converge in distribution in order for $\sqrt{\rho^\ast_n}\,M_n$ possessing a limit in distribution, see \cite[Thm.~2]{Topsoe}. As the sequence $\rho^\ast_n/\rho_n$ diverges, this is clearly not possible to hold.

Unfortunately an answer to the second question is not possible in this clarity, when considering non-trivial limits. Essentially, we can only say that the currently used methods do not allow for any conclusion on convergence. The limitations are the following: The central problem is that for the parameter range $\alpha\in[0,d]$ the current method does not provide tightness of the rescaled martingale sequence, hence we cannot infer that the sequence possesses a convergent subsequence. However, if tightness can be established in a different way then for the range $\alpha\in(\max\{1,d/2\},d]$ the limit has to be the diffusion process defined by the operator \eqref{def_covariance_op} as follows from the characterisation of any limit in the proof of the theorem. Here, the lower bound of $\max\{1,d/2\}$ results, on the one hand, from our estimation technique which necessitates $\alpha\geq 1$ and, on the other hand, from the definition of the limiting diffusion. Recall that the covariance operator is only of trace class for $\alpha>d/2$. Hence for $\alpha\in[0,d/2]$ we can no longer infer that the limiting diffusion even exists.
\end{remark}

\subsection{The mean-field Langevin equation}\label{section_langevin}

An important property of the limiting diffusion in view towards analytic and numerical studies is that it can be represented by a stochastic integral with respect to a cylindrical or $Q$--Wiener process. For a general discussion of infinite-dimensional stochastic integrals we refer to \cite{DaPratoZabczyk}. First, let $(W_t)_{t\geq 0}$ be a cylindrical Wiener process on $H^{-\alpha}(D)$ with covariance operator being the identity. Then, $G(\nu(t),t)\circ \iota^{-1}$ is a trace class operator on $H^{-\alpha}(D)$ for suitable values of $\alpha$. Here $\iota^{-1}:H^{-\alpha}(D)\to H^\alpha(D)$ is the Riesz Representation, i.e., the usual identification of a Hilbert space with its dual. The operator $G(\nu(t),t)\circ \iota^{-1}$ possesses a unique square-root we denote by $\sqrt{G(\nu(t),t)\circ\iota^{-1}}$ which is a Hilbert-Schmidt operator on $H^{-\alpha}(D)$. It follows that the stochastic integral process
\begin{equation}\label{stochastic_integral_process_1}
Z_t:=\int_0^t \sqrt{G(\nu(s),s)\circ\iota^{-1}}\,\dif W_s
\end{equation}
is a diffusion process in $H^{-\alpha}(D)$ with covariance operator $C(t)$. That is, $(Z_t)_{t\geq 0}$ is a version of the limiting diffusion in Theorem \ref{martingale_clt}. Now, formally substituting for the limits in \eqref{expansion_of_jump_process} yields the \emph{linear noise approximation}
\begin{equation*}
U_t=\nu_0+\int_0^t \tau^{-1}\,\bigl(U_s+\,F(U_s,s)\bigr)\,\dif s + \eps_n\int_0^t \sqrt{G(\nu(s),s)\circ\iota^{-1}}\,\dif W_s\,,
\end{equation*}
or in differential notation
\begin{equation}\label{linear_noise_approximation}
\dif U_t = \tau^{-1}\bigl(U_t+F(U_t,t)\bigr)\,\dif t +\eps_n\,\sqrt{G(\nu(t),t)\circ\iota^{-1}}\,\dif W_t\,,\quad U_0=\nu_0\,,
\end{equation}
where $\eps_n=\sqrt{v_+(n)/\ell_-(n)}$ is small for large $n$. Here we have used the operator notation
\begin{equation*}
F: H^{-\alpha}(D)\times\rr_+\to H^{-\alpha}(D):\,F(g,t)(x)= f\bigl(\langle g,w(x,\cdot)\rangle_{H^{\alpha}}+I(t,x)\bigr)\,.
\end{equation*}
Equation \eqref{linear_noise_approximation} is an infinite-dimensional stochastic differential equation with additive (linear) noise. Here additive means that the coefficient in the diffusion term does not depend on the solution $U_t$. A second formal substitution yields the \emph{Langevin approximation}. Here the dependence of the diffusion coefficient on the deterministic limit $\nu$ is formally substituted by a dependence on the solution. That is, we obtain a stochastic partial differential equation with multiplicative noise given by
\begin{equation*}
V_t=V_0+\int_0^t \tau^{-1}\,\bigl(V_s+\,F(V_s,s)\bigr)\,\dif s + \eps_n\int_0^t \sqrt{G(V_s,s)\circ\iota^{-1}}\,\dif W_s\,,
\end{equation*}
or in differential notation
\begin{equation}\label{langevin_approximation}
\dif V_t = \tau^{-1}\bigl(V_t+F(V_t,t)\bigr)\,\dif t +\eps_n\,\sqrt{G(V_t,t)\circ\iota^{-1}}\,\dif W_t\,.
\end{equation}
Note that the derivation of the above equations was only formal, hence we have to address the existence and uniqueness of solutions and the proper setting for these equations. This is left for future work. Furthermore, it is an ongoing discussion and probably undecidable as lacking a criterion of approximation quality which - if any at all - is the correct diffusion approximation to use. First of all note that for both versions the noise term vanishes for $n\to\infty$ and thus both have the Wilson-Cowan equation as their limit. And also, neither of them approximates even the first moment of the microscopic models exactly. This means that for neither we have that the mean solves the Wilson-Cowan equation which would be only the case if $f$ were linear. However, they are close to the mean of the discrete process. We discuss this aspect in the Appendix \ref{app_moment_equations}.

Furthermore, we already observe in the central limit theorem and thus also in the linear noise and Langevin approximation that the covariance \eqref{def_covariance_op} or the drift and the structure of the diffusion terms in \eqref{linear_noise_approximation} and \eqref{langevin_approximation}, respectively, are independent of objects resulting from the microscopic models. They are defined purely in terms of the macroscopic limit. This observation supports the conjecture that these approximations are independent from possible different microscopic models converging to the same deterministic limit. Analogous statements hold also for derivations from the van Kampen system size expansion \cite{Bressloff1} and in related limit theorems for reaction diffusion models \cite{Blount1,Kotelenez1,Kotelenez2}. 
The only object reminiscent of the microscopic models in the continuous approximations is the rescaling sequence $\eps_n$. However, the rescaling is proportional to the square root of $\ell_-(n)/v_+(n)$, i.e., the number of neurons per area divided by the size of the area, which is just the local density of particles. Therefore, in the approximations the noise scales inversely to the square root of neuron density in this model, which, interpreted in this way, can also be considered a macroscopic fixed parameter and chosen independently of the approximating sequence. 

\begin{remark}\label{remark_non-uniqueness_spde_approx} The stochastic partial differential equations \eqref{linear_noise_approximation} and  \eqref{langevin_approximation} which we proposed as the linear noise or Langevin approximation, respectively, are not necessarily unique as the representation of the limiting diffusion as a stochastic integral process \eqref{stochastic_integral_process_1} may not be unique. It will be subject for further research efforts to analyse the practical implications and usability of this Langevin approximation. Let $Q$ be a trace class operator,  $(W^Q_t)_{t\geq 0}$ be a $Q$-Wiener process and let $B(\nu(t),t)$ be operators such that $B(\nu(t),t)\circ Q\circ B(\nu(t),t)^\ast=G(\nu(t),t)\circ\iota^{-1}$, where $^\ast$ denotes the adjoint operator. Then also the stochastic integral process
\begin{equation*}
Z_t^Q:=\int_0^t B(\nu(s),s)\,\dif W_s^Q
\end{equation*}
is a version of the limiting diffusion in \eqref{martingale_clt} and the corresponding linear noise and Langevin approximations are given by
\begin{equation*}
\dif U_t^Q = \tau^{-1}\bigl(U_t^Q+F(U_t^Q,t)\bigr)\,\dif t +\eps_n\,B(\nu(t),t)\,\dif W_t^Q
\end{equation*}
and
\begin{equation*}
\dif V_t^Q = \tau^{-1}\bigl(V_t^Q+F(V_t^Q,t)\bigr)\,\dif t +\eps_n\,B(V^Q_t,t)\,\dif W_t^Q\,.
\end{equation*}
\end{remark}
\medskip

We conclude this section by presenting one particular choice of a diffusion coefficient and a Wiener process. We take $(W^Q_t)_{t\geq 0}$ to be a cylindrical Wiener process on $L^2(D)$ with covariance $Q=\textnormal{Id}_{L^2}$. Then we can chose $B(t)=j\circ(\cdot\sqrt{g(t)})\in L(L^2(D),H^{-\alpha}(D))$, where $j$ is the embedding operator $L^2(D)\hookrightarrow H^{-\alpha}(D)$ in the sense of \eqref{evolution_triple_prop}  and $(\cdot\sqrt{g(t)})\in L(L^2(D),L^2(D))$ denotes a pointwise product of a function in $L^2(D)$, i.e., 
\begin{equation*}
(\phi\cdot\sqrt{g(t)})(x)=\phi(x)\Bigl(\tau^{-1}\nu(t,x)+\tau^{-1}\,f\Bigl(\int_Dw(x,y)\,\nu(t,y)\,\dif y+I(t,x)\Bigr)\Bigr)^{1/2}\,.
\end{equation*}
We first investigate the operator $G(\nu(t),t)\circ\iota^{-1}$ and write it in more detail as the following composition of operators
\begin{equation*}
G(\nu(t),t)\circ\iota^{-1}=j\circ (\cdot g(t))\circ k \circ \iota^{-1}\,,
\end{equation*}
where $k$ is the embedding operator $H^\alpha(D)\hookrightarrow L^2(D)$. Next the Hilbert adjoint $B^\ast\in L(H^{-\alpha},L^2)$ is given by $B^\ast=(\cdot\sqrt{g})\circ k\circ\iota^{-1}$ which is easy to verify. Hence the stochastic integral of $B(t)$ with respect to $W^Q$ is again a version of the limiting martingale as
\begin{equation*}
B(t)\circ Q\circ B^\ast(t)=j\circ(\cdot\sqrt{g(t)})\circ \textnormal{Id}_{L^2}\circ(\cdot\sqrt{g(t)})\circ k\circ\iota^{-1}=j\circ (\cdot g(t))\circ k \circ \iota^{-1}=G(\nu(t),t)\circ\iota^{-1}\,.
\end{equation*}

\section{Discussion and extensions}\label{sec_conclusions}

In this article we have presented limit theorems that connect finite, discrete microscopic models of neural activity to the Wilson-Cowan neural field equation. The results state qualitative connections between the models formulated as precise probabilistic convergence concepts. Thus the results strengthen the connection derived in a heuristic way from the van Kampen system size expansion.

A general limitation of mathematically precise approaches to approximations, cf.~also the propagation to chaos limit theorems in \cite{Touboul},  is that the microscopic models are usually defined via the limit. In other words, the limit has to be known a-priori and we look for models which converge to this limit. Thus, in contrast to the van Kampen system size expansion the presented results are not a step-by-step modelling procedure in the sense that, via a constructive limiting procedure, a microscopic model yields a deterministic or stochastic approximation.
Hence, it might be objected that the presented method can only be used a-posteriori in order to justify a macroscopic model from a constructed microscopic model and that somehow one has to `guess' the correct limit in advance. Several remarks can be made to answer this objection.

%
%
%
First, this observation is certainly true, but not necessarily a drawback. On the contrary; when both microscopic and macroscopic models are available, then it is rather important to know how these are connected and qualitatively and quantitatively characterise this connection. 
Concerning neural field models, this precise connection was simply not available so far for the well-established Wilson-Cowan model. Furthermore, when starting from a stochastic microscopic description working through proving the conditions for convergence for given microscopic models one obtains very strong hints on the structure of a possible deterministic limit. Therefore our results can also ease the procedure of `guessing the correct limit'.

%
%
%
%
Secondly, often a phenomenological, deterministic model which is an approximation to an inherently probabilistic process is derived from ad-hoc heuristic arguments. Given that the model has proved useful, one often aims to derive a justification from first principles and / or a stochastic version which keeps the features of the deterministic model but also accounts for the formerly neglected fluctuations. 
A standard, though somewhat simple approach to obtain stochastic versions consists of adding (small) noise to the deterministic equations. This article, provides a second approach which consists of finding microscopic models, which converge to the deterministic limit to obtain a stochastic correction via a central limit argument.

%
Thirdly and finally, the method also provides an argument for new equations, i.e., the Langevin and linear noise approximations, which can be used to study the stochastic fluctuations in the model. Furthermore, in contrast to previous studies we do not provide deterministic moment equations but stochastic processes, which can be, e.g., via Monte Carlo simulations, studied concerning a large number of pathwise properties and dynamics beyond first and second moments.

\medskip

We now conclude this article commenting on the feasibility of our approach connecting microscopic Markov models to deterministic macroscopic equations when dealing with different master equation formulations that appear in the literature. Additionally, the following discussions also relate the model \eqref{master_equation_2} considered in this article to other master equation formulations. We conjecture that the analogous results as presented for the Wilson-Cowan equation \eqref{Wilson_Cowan} in Section \ref{sec_precise_LLN} also hold for these variations of the master equations. This should be possible to achieve by an adaptation of the methods of proof presented although we have not performed the computations in detail.

\subsection{A Variation of the master equation formulation}

A first variation of the discrete model we discussed in Section \ref{sec_master_equations} was considered in the articles \cite{BuiceCowan,BuiceCowanChow} and a version restricted to a bounded state space also appears in \cite{TouboulErmentrout}. This model consists of the master equation stated below in \eqref{master_equation_effective_spike} which closely resembles \eqref{master_equation_2}. In the earlier reference \cite{BuiceCowan} the model was introduced with a different interpretation called the \emph{effective spike model}. We briefly explain this interpretation before presenting the master equation. Instead of interpreting $P$ as the number of neuron populations,  in this model $P$ denotes the number of different neurons in the network located within a spatial domain $D$. Then $\Theta^k_t$, the state of the $k$th neuron, counts the number of `effective' spikes this neuron has emitted in the past up till time $t$. Effective spikes are those spikes that still influence the dynamics of the system, e.g., via a post-synaptic potential. Then state transitions adding / subtracting one effective spike for the $k$th neuron are governed by a firing rate function $\widetilde f_k$, which depends on the input into neuron $k$, and a decay rate $\tau^{-1}$. The constant decay rate indicates that emitted spikes are effective for a time interval of length $\tau$ and the gain function is defined -- neglecting external input -- by
\begin{equation*}
\widetilde f_k(\theta)=f^\ast\Bigl(\sum_{j=1}^P \widetilde W_{kj}\theta^j\Bigr)\,, 
\end{equation*}
where $f^\ast$ is a certain nonnegative, real function. It is stated clearly in \cite{BuiceCowanChow} that the function $f^\ast$ is \emph{not} equal to the gain function $f$ in the proposed limiting Wilson-Cowan equation \eqref{Wilson_Cowan} but rather connected to $f$ such that
\begin{equation}\label{conjectured_connection_0}
\EX f^\ast\Bigl(\sum_{j=1}^P \widetilde W_{kj}\Theta^j_t\Bigr)=f\Bigl(\sum_{j=1}^P \widetilde W_{kj}\EX\Theta^j_t\Bigr) +\textnormal{ higher order terms}\,.
\end{equation}
The authors in \cite{BuiceCowanChow} state that for any function $f$ such a function $f^\ast$ can be found. 
Then the process $\Theta_t=(\Theta^1_t,\ldots,\Theta^P_t)$ is a jump Markov process with its evolution governed by the master equation
\begin{equation}\label{master_equation_effective_spike}
\frac{\dif\Pr[\theta,t]}{\dif t}=\sum_{k=1}^P \Bigl[\widetilde f_k(\theta-e_k)\,\Pr[\theta-e_k,t]-\Bigl(\frac{1}{\tau}\,\theta^k+\widetilde f_k(\theta)\Bigr)\,\Pr[\theta,t]+\frac{1}{\tau}\,(\theta^k+1)\,\Pr[\theta+e_k,t]\Bigr]
\end{equation}
with boundary conditions $\Pr[\theta,t]=0$ if $\theta\notin \mathbb{N}_0^P$ as stated in \cite{BuiceCowanChow}. The advantage of the effective spike model interpretation over the interpretation as neurons per population is that the unbounded state space of the model is justified. In principle there can be an arbitrary number of spikes emitted in the past still active. However, a disadvantage of the master equation \eqref{master_equation_effective_spike} is that for taking the limit it lacks a parameter corresponding to the system size providing a natural small parameter in the van Kampen system size expansion. This explains the shift in the interpretation of the master equation in the study \cite{BuiceCowanChow} following \cite{BuiceCowan} and subsequently in \cite{Bressloff1} to the interpretation we presented in Section \ref{sec_master_equations} which provides the system-size parameters $l(k)$.

On the level of Markov jump processes the master equation \eqref{master_equation_effective_spike} obviously describes dynamics similar to the master equation \eqref{master_equation_2} only replacing the activation rate $\tau^{-1} l(k)\nbar f_k(\theta)$ in \eqref{master_equation_2} by $\widetilde f_k(\theta)$ which is independent of the parameter $l(k)$. Thus, the model \eqref{master_equation_effective_spike} can be understood as resulting from \eqref{master_equation_2} \emph{after} a limit procedure taking $l(k)\to\infty$ has been applied and the firing rate functions are connected via the formal limit $\lim_{l(k)\to\infty} l(k)\nbar f_k(\theta) =\widetilde f_k(\theta)$. A qualitative interpretation of this limit procedure connecting the two types of models is given in \cite{BuiceCowan}. This observation motivated the model in \cite{Bressloff1} stepping back one limit procedure and thus providing the correct framework for the derivation of limit theorems.

It would be an interesting addition to the limit theorems in Theorem \ref{lln} to derive a law of large numbers for the models \eqref{master_equation_effective_spike} with stochastic mean activity $\nu^n$ as defined in \eqref{def_coordinate_function} and suitable chosen weights $\ntilde W_{kj}$. Clearly, the macroscopic limit should be given by the Wilson-Cowan equation \eqref{Wilson_Cowan}. We conjecture that the appropriate condition for the function $f^\ast$ in the present setting -- including time dependent inputs -- is
\begin{equation}\label{conjectured_connection}
\EX\Bigl[ l(k,n)^{-1} f^\ast\Bigl(\sum_{j=1}^P \widetilde W_{kj}^n\Theta^j+\widetilde I_{k,n}(t)\Bigr)\Bigr]=f\Bigl(\sum_{j=1}^P \nbar W_{kj}^n\frac{\EX\Theta_j}{l(j,n)}+\nbar I_{k,n}(t)\Bigr) +\textnormal{ h.o.t.}\,,
\end{equation}
such that the higher order terms are uniformly bounded and vanish in the limit $n\to\infty$, and where the weights $\nbar W_{kj}^n$ and inputs $\nbar I_{k,n}(t)$ are defined as in \eqref{definition_of_discret_weights} and \eqref{definition_of_discret_inputs}. Property \eqref{conjectured_connection} closely resembles condition \eqref{conjectured_connection_0} and trivially holds for linear $f$ with $f^\ast=f$.

\subsection{Bounded state space master equations}\label{section_discussion_bounded_state_space}

We have already stated when introducing the microscopic model in Section \ref{sec_master_equations} that the interpretation of the parameter $l(k)$ as the number of neurons in the $k$-th population is not literally correct. The state space of the process is unbounded, hence arbitrarily many neurons can be active and thus each population contains arbitrarily many neurons. In order to overcome this interpretation problem it was supposed to consider the master equation only on a bounded state space. That is, the $k$-th population consists of $l(k)$ neurons and $0\leq \Theta^k_t\leq l(k)$ almost surely. Such master equations are simply obtained by setting the transition rates for transition of $\theta^k$ from $l(k)\to l(k)+1$ to zero.

A first master equation of this form was considered in \cite{OhiraCowan2} which, in present notation, takes the form
\begin{equation}
\frac{\dif\Pr[\theta,t]}{\dif t}=\frac{1}{\tau}\sum_{k=1}^P\Bigl[(l(k)-\theta^k+1)\nbar f_k(\theta-e_k)\Pr[\theta^k-e_k,t] -\Bigl(\theta^k+(l(k)-\theta^k)\nbar f_k(\theta)\Bigr)\Pr[\theta,t] + (\theta^k+1)\Pr[\theta+e_k,t] \Bigr]\,.\label{master_equation_3}
\end{equation}
Versions of such a master equation for, e.g., one population only or coupled inhibitory and excitatory populations were considered in \cite{Benayoun,OhiraCowan2} and a van Kampen systems size expansion was carried out. Here the bound in the state space provides a natural parameter for the rescaling, thus a small parameter for the expansion. The setup of this problem resembles closely the structure of excitable membranes for which limits have been obtained with the present technique by one of the present author and co-workers in \cite{RTW}. Therefore we conjecture that our limit theorems also apply to this setting with minor adaptations with essentially the same conditions and results as in Section \ref{sec_precise_LLN}. However, the macroscopic limit which will be obtained does not conform with the Wilson-Cowan equation but will be given by
\begin{equation}\label{really_rate_equation}
\tau\,\dot\nu(t,x)=-\nu(t,x)+(1-\nu(t,x))f\Bigl(\int_D w(x,y)\nu(t,y)\,\dif y +I(t,x)\Bigr)\,.
\end{equation}

Next, we return to the master equation \eqref{master_equation_2} as discussed in this article in Section \ref{sec_master_equations} and the comment we made regarding bounded state spaces the footnote on page \pageref{the_bounded_state_space_footnote}. In our primary reference for this model \cite{Bressloff1} actually a bounded state space version of the master equation was considered where the activation rate for the event $\theta^k\to\theta^k+1$ is
\begin{equation}\label{bressloff_discontinuous_rate}
l(k)\nbar f_k(\theta,t)\mathbb{I}_{[\theta_k<l(k)]}\,, 
\end{equation}
replacing $l(k)\nbar f_k(\theta,t)$ in \eqref{master_equation_2}. The van Kampen system size expansion was then applied to this bounded state space master equation, tacitly neglecting possible difficulties which might arise due to the discontinuity of \eqref{bressloff_discontinuous_rate} considered as a function on $\rr^P$. However, for the present, mathematically precise limit convergence results considering bounded state space as originally suggested in \cite{Bressloff1} are problematic. 
The discontinuous activation rate \eqref{bressloff_discontinuous_rate} causes the machinery developed in \cite{RTW} which depends on Lipschitz-type estimates to break down. 
However, we strongly expect that also in this case the law of large numbers with the deterministic limit given by the Wilson-Cowan equation \eqref{Wilson_Cowan} holds. Furthermore, also the Langevin approximations should agree with the equations discussed in Section \ref{section_langevin}. However, we have not yet been able to prove such a theorem. We further conjecture that the results in this article can be used to prove the convergence for the bounded state space model by a domination argument. Heuristically, it seems clear that a bounded process should be dominated by a process that possesses the same dynamics inside the state space of the bounded process but can stray out from that bounded domain. Hence, as the limit of the potentially larger process lies within the domain where the two processes agree also the dominated process should converge to the same limit. Mathematically, this line of argument relies on non-trivial estimates between occupation measures of high-dimensional Markov processes. This is work in progress.

\subsection{Activity based neural field model}

Finally, we return also to a difference in neural field theory mentioned in the beginning. In contrast to rate-based neural field models of the Wilson-Cowan type \eqref{intro_wc_eq} there exists a second essential class of neural field models, so called activity based models, the prototype of which is the \emph{Amari equation}
\begin{equation}\label{Amari_equation}
\tau\,\dot\nu(t,x)=-\nu(t,x)+\int_D w(x,y)f\bigr(\nu(t,y)\bigr)\,\dif y+I(t,x)\,.
\end{equation}
We conjecture that also for this type of equations a phenomenological microscopic model can be constructed with a suitable adaptation of the activation rates and that limit theorems analogous to the results in Section \ref{section_llns} hold. Then also a Langevin equation for this model can be obtained and used for further analysis.

\section{Proofs of the main results}\label{section_proofs}

In this section we present the proofs of the limit theorems. For the convenience of the reader, as it is important tool in the subsequent proofs, we first state the Poincar\'e inequality. Let $D\subset\rr^d$ be a convex domain then it holds for any function $\phi\in H^1(D)$ that
\begin{equation}\label{poincare_inequality} 
\|\nbar\phi_D-\phi\|_{L^2(D)}\leq \frac{\textnormal{diam}(D)}{\pi}\,\|\nabla\phi\|_{L^2(D)}\,,
\end{equation}
where $\nbar \phi_D$ is the mean value of the function $\phi$ on the domain $D$, i.e.,
\begin{equation}\label{overall_proof_pwc_approx}
\nbar \phi_D=\frac{1}{|D|}\int_D \phi(x)\,\dif x\,.
\end{equation}
Moreover, the constant in the right hand side of \eqref{poincare_inequality} is the optimal constant depending only on the diameter of the domain $D$, cf.~\cite{Acosta,Payne}. Whenever we omit to denote the spatial domain for definition of norms or inner products in $L^2(D)$ or Sobolev spaces $H^\alpha(D)$ then it is to be interpreted as the norm over the whole domain $D$. If the norm is taken only over a subset $D_{k,n}$ then this is always indicated unexceptionally.
\medskip

For the benefit of the reader we next repeat the limiting equation
\begin{equation}\label{Wilson_Cowan_2}
\tau\,\dot\nu(t,x)=-\nu(t,x)+f\Bigl(\int_D w(x,y)\nu(t,y)\,\dif y+I(t,x)\Bigr)\,.
\end{equation}
We denote by $F$ the Nemytzkii operator on $L^2(D)$ defined by
\begin{equation}\label{def_Nemytzkii_op}
F(g,t)(x)= f\Bigl(\int_D w(x,y)g(y)\,\dif y+I(t,x)\Bigr)\qquad\forall\,g\in L^2(D)\,,
\end{equation}
and for all $\theta\in\mathbb{N}^P_0$ we define a discrete version of the Nemyztkii operator via
\begin{eqnarray}
-\frac{1}{\tau}\,\nu^n(\theta)+\frac{1}{\tau}\,\nbar F^n(\nu^n(\theta),t)&=&\lambda^n(\theta,t)\int_{\mathbb{N}_0^P} \bigl(\nu^n(\xi)-\nu^n(\theta)\bigr)\,\mu^{n}\bigl((\theta,t),\dif\xi\bigr)\nonumber\\
&=& \frac{1}{\tau}\,\sum_{k=1}^{P} \frac{1}{l(k,n)}\Bigl(-\theta^{k}+l(k,n)\,\nbar f_{k,n}(\theta,t)\Bigr)\,\mathbb{I}_{D_{k,n}}\nonumber\\
&=& -\frac{1}{\tau}\,\nu^n(\theta)+\frac{1}{\tau} \sum_{k=1}^{P} \nbar f_{k,n}(\theta,t)\,\mathbb{I}_{D_{k,n}}\,.\phantom{xxx}\label{def_discrete_nemytzkii}
\end{eqnarray}
Note that $\tau^{-1}(\phi,\nu^n(\theta))_{L^2}+\tau^{-1}(\phi,\nbar F^n(\nu^n(\theta),t))_{L^2}$ for $\phi\in L^2(D)$ corresponds to the generator of $(\Theta^n_t,t)_{t\geq0}$ applied to the function $(\theta,t)\mapsto (\phi,\nu^n(\theta))_{L^2}$.

\medskip

Finally, another useful property is that the means of the process' components are bounded. For each $k,n$ it holds that
\begin{equation*}
 \EX \Theta^{k,n}_t\,=\,\EX \Theta^{k,n}_0+\frac{1}{\tau}\int_0^t l(k,n)\,\EX \nbar f_{k,n}(Y^n_s)-\EX\Theta^{k,n}_s\,\dif s\,\leq\, \,\EX \Theta^{k,n}_0+\frac{1}{\tau}\int_0^t l(k,n)\,\|f\|_0-\EX\Theta^{k,n}_s\,\dif s,
\end{equation*}
see also \eqref{mean_eq_MC}. Therefore it holds that $\EX\Theta^{k,n}_t\leq m_t^{k,n}$, where $m_t^{k,n}$ solves the deterministic initial value problem
\begin{equation*}
 \dot m^{k,n}_t = -\frac{1}{\tau}\, m^{k,n}_t + \frac{1}{\tau}\, l(k,n)\,\|f\|_0,\quad m_0^{k,n}=\EX\Theta^{k,n}_0,
\end{equation*}
i.e., 
\begin{equation}\label{components_mean_bound}
m_t^{k,n}=\e^{-t/\tau}\bigl(m_{k,n}^0-l(k,n)\|f\|_0\bigl)+l(k,n)\|f\|_0\,\leq\,l(k,n)\bigl(1+\|f\|_0\bigr) \quad\forall\,t\geq 0\,.
\end{equation}
Here we also used the assumption $\EX^n\Theta^{k,n}_0\leq l(k,n)$ on the initial condition.

\subsection{Proof of Theorem \ref{lln} (Law of large numbers)}\label{section_proof_lln}

In order to prove the law of large numbers, Theorem \ref{lln}, we apply the law of large numbers for Hilbert space valued PDMPs, see \cite[Thm.~4.1]{RTW}, to the sequence of homogeneous PDMPs $(Y_t^n)_{t\geq 0}=(\Theta^n_t,t)_{t\geq 0}$. For the application of this theorem, recall that the first, piecewise constant, vector-valued component of this process counts the number of active neurons in each sub-population and the second, deterministic component states time. The process $(Y_t^n)_{t\geq 0}$ is the usual `space-time process', i.e., homogeneous Markov process which is obtained via a state-space extension to obtain a homogeneous Markov process from the inhomogeneous process $(\Theta^n_t)_{t\geq0}$. The continuous component satisfies the simple ODE $\dot t=1$, $t(0)=0$ and thus the full process is a PDMP. In the terminology of \cite{RTW} the sequence of coordinate functions on the different state spaces of the PDMPs $(Y_t^n)_{t\geq 0}$ into a common Hilbert space is given by the maps $\nu^n$ \eqref{def_coordinate_function} with the common Hilbert space $L^2(D)$. Thus in order to infer convergence in probability \eqref{thm_lln_conv_in_prob} from \cite[Thm.~4.1]{RTW} it is sufficient to validate the following conditions: 
\setlength{\leftmargini}{4.5em}
\begin{itemize}[align=left]
  \setlength{\labelwidth}{4.0em}
  \setlength{\labelsep}{0.5em}
\item[\textnormal{(LLN1)}] For fixed $T>0$ it holds that
\begin{equation}\label{proof_lln_cond_a}
\lim_{n\to\infty}\EX^n\int_0^T\lambda^n(Y^n_t)\int_{\mathbb{N}^P} \|\nu^n(\xi)-\nu^n(\Theta^n_t)\|_{L^2(D)}^2\,\mu^n(Y_t^n,\dif\xi)\,\dif t=0\,.
\end{equation}
\item[\textnormal{(LLN2)}] The Nemytzkii operator $F$ satisfies a Lipschitz condition in $L^2(D)$ uniformly with respect to $t$, $t\geq 0$, i.e., there exists a constant $L_0>0$ such that
\begin{equation}\label{lipschitz_condition_in_H0}
\|F(g_1,t)-F(g_2,t)\|_{L^2}\,\leq \,L_0\,\|g_1-g_2\|_{L^2}\qquad \forall\, t\geq 0,\, g_1,g_2\in L^2(D)\,.
\end{equation}
\item[\textnormal{(LLN3)}] For fixed $T>0$ it holds that
\begin{equation}\label{finite_var_bound}
\lim_{n\to\infty}\EX^n\int_0^T\big\|\nbar F^n(\nu^n_t,t) -  F(\nu^n_t,t)\big\|_{L^2}\,\dif t = 0\,.
\end{equation}
\end{itemize}
Note that the final condition of \cite[Thm.~4.1]{RTW}, i.e., the convergence of the initial conditions, is satisfied by assumption. For a discussion of these conditions we refer to \cite{RTW} and proceed to their derivation for the present model in the subsequent parts (a) to (c).\medskip

(a)\quad In order to prove condition \eqref{proof_lln_cond_a} we write the integral with respect to the discrete probability measure $\mu^n$ as a sum. This yields
\begin{eqnarray}\label{lln_martingale_estimate}
\lefteqn{\EX^n\lambda(Y^n_t)\int_{\mathbb{N}^P} \|\nu^n(\xi)-\nu^n(\Theta^n_t)\|_{L^2(D)}^2\,\mu^n(Y_t^n,\dif\xi)}\nonumber\\
&&\phantom{xxxxxxxxxxxx}=\ \frac{1}{\tau}\sum_{k=1}^P\EX^n \frac{1}{l(k,n)^2}\Bigl(\Theta^{k,n}_t+l(k,n)\,\nbar f_{k,n}(Y^n_t)\Bigr)\, |D_{k,n}|\\
&&\phantom{xxxxxxxxxxxx}\leq\ \frac{1}{\tau}\,\frac{1 +2\|f\|_0}{\ell_-(n)}\,|D|\,,\nonumber
\end{eqnarray}
where we have used the upper bound \eqref{components_mean_bound} on the expectation $\EX^n \Theta^{k,n}_t$ and the assumption on the initial conditions. Next, integrating over $[0,T]$ and employing the assumption $\lim_{n\to\infty}\ell_-(n)=\infty$ in \eqref{assumptions_lln} establishes condition \eqref{proof_lln_cond_a}.\medskip

(b)\quad The Lipschitz condition \eqref{lipschitz_condition_in_H0} of the Nemytzkii operators is a straightforward consequence of the Lipschitz continuity \eqref{lipschitz_on_f} of the gain function $f$ as
\begin{eqnarray*}
\|F(g_1,t)-F(g_2,t)\|_{L^2}^2&=&\int_D \Big|f\Bigl(\int_D w(x,y)g_1(y)\dif y+I(x,t)\Bigr)-f\Bigl(\int_D w(x,y)g_2(y)\dif y+I(x,t)\Bigr)\Big|^2\,\dif x\\
&\leq&  L^2\,\int_D\Big|\int_D w(x,y)\bigr(g_1(y)-)g_2(y)\bigl)\dif y\Big|^2\,\dif x\\
&\leq& L^2\,\int_D \|w(x,\cdot)\|_{L^2}^2\,\|g_1-g_2\|_{L^2}^2\,\dif x\\
&=& L^2\,\|w\|_{L^2\times L^2}^2\,\|g_1-g_2\|_{L^2}^2\,.
\end{eqnarray*}
Therefore \eqref{lipschitz_condition_in_H0} holds with Lipschitz constant $L_0:=L\,\|w\|_{L^2\times L^2}$.\medskip

(c)\quad Finally we prove the convergence of the generators \eqref{finite_var_bound}. To this end we employ the characterisation of the norm in $L^2(D)$ by $\|\eta\|_{L^2}=\sup_{\|\phi\|_{L^2}=1}\big|(\phi,\eta)_{L^2}\big|$ for all $\eta\in L^2(D)$ and thus consider first the scalar product of elements $\phi\in L^2(D)$ with $\|\phi\|_{L^2}=1$ and the difference inside the norm in \eqref{finite_var_bound}. On the one hand we obtain using definition \eqref{def_discrete_nemytzkii} that
\begin{eqnarray}
\bigl(\phi,\nbar F^n(\nu^n_t,t)\bigl)_{L^2}&=&\Bigl(\phi,\sum_{k=1}^P \nbar f_{k,n}(Y^n_t)\,\mathbb{I}_{D_{k,n}}\Bigr)_{L^2}\,.\phantom{xxx}\label{difference_part_1}
\end{eqnarray}
Next we apply the Nemytzkii operator $F$ defined in \eqref{def_Nemytzkii_op} to $\nu^n(t)$ and take the inner product of the result with respect to $\phi$ to obtain on the other hand
\begin{equation}
\bigl(\phi,F(\nu^n_t,t)\bigr)_{L^2}=\Bigl(\phi,f\Bigl(\sum_{k=1}^P \frac{\Theta^{k,n}_t}{l(k,n)}\int_{D_{k,n}} w(\cdot,y)\,\dif y +I(t,\cdot)\Bigr)\,.\label{difference_part_2}
\end{equation}

Subtracting \eqref{difference_part_2} from \eqref{difference_part_1} we obtain the integrated difference
\begin{eqnarray*}
\lefteqn{\hspace{-10pt}\bigl(\phi,\nbar F^n(\nu^n_t,t)\bigl)_{L^2}-\bigl(\phi,F(\nu^n_t,t)\bigr)_{L^2}\ =}\\[2ex]
&=&\int_D \phi(x)\,\biggl[
\sum_{k=1}^P \nbar f_{k,n}(Y^n_t)\,\mathbb{I}_{D_{k,n}}-
f\Bigl(\sum_{j=1}^P \frac{\Theta^n_j(t)}{l(j,n)}\int_{D_{j,n}} w(x,y)\,\dif y +I(t,x)\Bigr)\biggr]
\,\dif x\\
&=&\sum_{k=1}^P\int_{D_{k,n}} \phi(x)\,\biggl[
f\Bigl(\sum_{j=1}^P \nbar W_{kj}^n\frac{\Theta^{j,n}_t}{l(j,n)} + \nbar I_{k,n}(t)\Bigr)-
f\Bigl(\sum_{j=1}^P \frac{\Theta^{j,n}_t}{l(j,n)}\int_{D_{j,n}} w(x,y)\,\dif y +I(t,x)\Bigr)\biggr]\,\dif x\,.
\end{eqnarray*}
We proceed to estimate the norm of the term in the right hand side. We use the Lipschitz condition \eqref{lipschitz_on_f} on $f$, the triangle inequality and finally the Cauchy-Schwarz inequality on the resulting second term to obtain the estimate
\begin{eqnarray*}
 \lefteqn{\big|\bigl(\phi,\nbar F^n(\nu^n_t,t)\bigl)_{L^2}-\bigl(\phi,F(\nu^n_t,t)\bigr)_{L^2}\big|}\\[2ex]
&\leq&L\sum_{k=1}^P\int_{D_{k,n}}|\phi(x)|\,\Big|\sum_{j=1}^P\frac{\Theta^n_j(t)}{l(j,n)} \Bigl(\nbar W_{kj}^{n}-\int_{D_{j,n}} w(x,y)\,\dif y\Bigr) + \nbar I_{k,n}(t) - I(t,x) \Big|\,\dif x\\
&\leq & L\underbrace{\sum_{k=1}^P\int_{D_{k,n}}|\phi(x)|\,\Big|\sum_{j=1}^P\frac{\Theta^n_j(t)}{l(j,n)} \Bigl(\nbar W_{kj}^{n}-\int_{D_{j,n}} w(x,y)\,\dif y\Bigr) \Big|\,\dif x}_{(\ast)} + L\underbrace{\sum_{k=1}^P \|\phi\|_{L^2(D_{k,n})}\,\|\nbar I_{k,n}(t) - I(t)\|_{L^2(D_{k,n})}}_{(\ast\ast)}.
\end{eqnarray*}
Here, the term in the right hand side marked $(\ast\ast)$ is further estimated using the Cauchy-Schwarz inequality and the Poincar\'e inequality \eqref{poincare_inequality} which yields
\begin{eqnarray}\label{intermediate_estimate_1}
(\ast\ast)\ \leq \ \frac{\delta_+(n)}{\pi}\Bigl(\sum_{k=1}^P\|\nabla I(t)\|_{L^2(D_{k,n})}^2 \Bigr)^{1/2}\ = \ \frac{\delta_+(n)}{\pi} \|\nabla_{\!x} I(t)\|_{L^2}\,.
\end{eqnarray}
We now consider the term marked $(\ast)$. Inserting the definition of $\nbar W^{n}_{kj}$ given in \eqref{definition_of_discret_weights},  the reordering of the summations and changing the order of integration yields
\begin{eqnarray*}
(\ast)&=& \sum_{k=1}^P\int_{D_{k,n}}|\phi(x)|\,\bigg|\sum_{j=1}^P\frac{\Theta^{j,n}_t}{l(k,n)} \biggl(\int_{D_{j,n}}\Bigl(\frac{1}{|D_{k,n}|}\int_{D_{k,n}}w(z,y)\,\dif z\Bigr)- w(x,y)\,\dif y\biggr) \bigg|\,\dif x\\[1ex]
&\leq&\sum_{k=1}^P\int_{D_{k,n}}|\phi(x)|\,\sum_{j=1}^P\frac{\Theta^{j,n}_t}{l(k,n)} \int_{D_{j,n}}\Big|\Bigl(\frac{1}{|D_{k,n}|}\int_{D_{k,n}}w(z,y)\,\dif z\Bigr)- w(x,y)\Big|\,\dif y\,\dif x\\[1ex]
&=& \sum_{k=1}^P\sum_{j=1}^P\frac{\Theta^{j,n}_t}{l(k,n)}\int_{D_{j,n}}\Biggl[\int_{D_{k,n}}|\phi(x)|\,\Big|\Bigl(\frac{1}{|D_{k,n}|}\int_{D_{k,n}}w(z,y)\,\dif z\Bigr)- w(x,y)\Big|\,\dif x\Biggr]\,\dif y\,.
\end{eqnarray*}
We next apply the Cauchy-Schwarz inequality to the integral inside the square brackets in the last term. Thus we obtain the estimate
\begin{eqnarray*}
(\ast)&\leq& \sum_{j=1}^P\int_{D_{j,n}}\sum_{k=1}^P\|\phi\|_{L^2(D_{k,n})}\biggl[\int_{D_{k,n}}\,\Big|\Bigl(\frac{1}{|D_{k,n}|}\int_{D_{k,n}}w(z,y)\,\dif z\Bigr)- w(x,y)\Big|^2\,\dif x\biggr]^{1/2}\,\dif y\,.
\end{eqnarray*}
Now the Poincar\'e inequality \eqref{poincare_inequality} is applied to the innermost integral inside the square brackets which yields
\begin{eqnarray*}
(\ast)&\leq& \sum_{j=1}^P \frac{\Theta^{j,n}_t}{l(k,n)}\int_{D_{j,n}}\sum_{k=1}^P\|\phi\|_{L^2(D_{k,n})}\, \frac{\textnormal{diam}(D_{k,n})}{\pi}\,\|\nabla_{\!x} w(\cdot,y)\|_{L^2(D_{k,n})}\,\dif y\,.
\end{eqnarray*}
Finally, using once more the Cauchy-Schwarz inequality on the innermost summation we obtain
\begin{eqnarray}\label{intermediate_estimate_2}
(\ast)&\leq& \frac{\delta_+(n)}{\pi}\sum_{j=1}^P \frac{\Theta^{j,n}_t}{l(k,n)}\int_{D_{j,n}} \|\nabla_{\!x} w(\cdot,y)\|_{L^2}\,\dif y\,.
\end{eqnarray}

Now, a combination of the estimates \eqref{intermediate_estimate_1} and \eqref{intermediate_estimate_2} on the terms $(\ast)$ and $(\ast\ast)$ yields
\begin{equation*}
\Bigl|\bigl(\phi,\nbar F^n(\nu^n_t,t)\bigr)_{L^2}-\bigl(\phi,F(\nu^n(t),t)\bigr)_{L^2}\Bigr|\,\leq\, \delta_+(n)\,\frac{L}{\pi}\biggl(\sum_{j=1}^P \frac{\Theta^{j,n}_t}{l(k,n)}\int_{D_{j,n}} \|\nabla_{\!x} w(\cdot,y)\|_{L^2}\,\dif y+\|\nabla_{\!x} I(t)\|_{L^2}\biggr)\,.
\end{equation*}
Here the right hand side is independent of $\phi$, hence taking the supremum over all $\phi$ with $\|\phi\|_{L^2}=1$ yields
\begin{equation*}
\bigl\|\nbar F^n(\nu^n_t,t)-F(\nu^n(t),t)\bigr\|_{L^2}\,\leq\, \delta_+(n)\,\frac{L}{\pi}\biggl(\sum_{j=1}^P \frac{\Theta^{j,n}_t}{l(k,n)}\int_{D_{j,n}} \|\nabla_{\!x} w(\cdot,y)\|_{L^2}\,\dif y+\|\nabla_{\!x} I(t)\|_{L^2}\biggr)\,.
\end{equation*}
Finally, integrating over $(0,T)$ and taking the expectation on both sides results in
\begin{eqnarray}\label{final_estimate_in_finite_var_proof}
\EX^n\!\!\int_0^T\bigl\|\nbar F^n(\nu^n_t,t)-F(\nu^n(t),t)\bigr\|_{L^2}\dif t&\leq& \delta_+(n)\,\frac{L}{\pi}\Bigl(\sqrt{|D|}\,T\,\bigl(1+\|f\|_0\bigr)\,\|\nabla_{\!x} w\|_{L^2\times L^2}+\|\nabla_{\!x} I\|_{L^1((0,T),L^2)}\Bigr)\,.\nonumber\\
\end{eqnarray}
Here we have used \eqref{components_mean_bound} and a combination of the Cauchy-Schwarz and Poincar\'e inequality \eqref{poincare_inequality} in order to estimate
\begin{equation*}
\EX^n \sum_{j=1}^P \frac{\Theta^{j,n}_t}{l(k,n)}\int_{D_{j,n}} \|\nabla_{\!x} w(\cdot,y)\|_{L^2}\,\dif y\,\leq \, \delta_+(n)\,\frac{\sqrt{|D|}\,\bigl(1+\|f\|_0\bigr)}{\pi}\, \|\nabla_{\!x} w\|_{L^2\times L^2}\,.
\end{equation*}
The upper bound in \eqref{final_estimate_in_finite_var_proof} is of order $\landau(\delta_+(n))$ and therefore converges to zero for $n\to\infty$ due to assumption \eqref{assumptions_lln}. Hence, condition \eqref{finite_var_bound} is satisfied as convergence in the mean implies convergence in probability. The proof of the convergence in probility \eqref{thm_lln_conv_in_prob} is completed.\medskip

It is now easy to extend this result to the convergence in the $r$-th mean. First of all the convergence in probility \eqref{thm_lln_conv_in_prob} implies for all $r\geq 1$ the convergence in probability of the random variables $\sup_{t\in[0,T]}\|\nu^n_t-\nu(t)\|_{L^2}^r$ to zero. As convergence in the mean of real valued random variables is equivalent to convergence in probability and uniform integrability it remains to prove the latter for the families $\sup_{t\in[0,T]}\|\nu^n_t-\nu(t)\|_{L^2}^r$, $n\in\mathbb{N}$. 

We first consider the case $r=1$, and establish a uniform bound on the second moments $\EX^n\sup_{t\in[0,T]}\|\nu^n_t-\nu(t)\|_{L^2}^2$. Then the de la Vall\'ee-Poussin Theorem, cf.~\cite[App., Prop.~2.2]{EthierKurtz}, implies that the random variables $\sup\nolimits_{t\in[0,T]}\|\nu^n_t-\nu(t)\|_{L^2}$, $n\in\mathbb{N}$, are uniformly integrable.

Without loss of generality we can assume that there exist\footnote{The Poisson process jumps at a faster rate than the components of the Markov chain regardless of the time and the state these are in. Furthermore all jumps are upwards. Hence using a coupling argument as discussed in the proof of \cite[Thm.~4.3.5]{Jacobsen} we find that there exists a probability space supporting two processes with distributions equivalent to the Poisson process and the Markov chain component such that the Poisson process dominates the second process for all paths. Clearly, all moments dominate and this inequalities are then valid for any probability spaces supporting these processes.}  Poisson processes $(N^{k,n}_t)_{t\geq 0}$ with rates $\Lambda_{k,n}=l(k,n)(1+\|f\|_0)/\tau$, which dominate $(\Theta^{k,n}_t-\Theta^{k,n}_0)_{t\geq 0}$ pathwise. Then we obtain almost surely
\begin{eqnarray*}
\|\nu^n_t\|_{L^2}^2\,\leq\, 2\|\nu^n_0\|_{L^2}^2+2\sum_{k=1}^P\frac{(\Theta^{k,n}_t-\Theta^{k,n}_0)^2}{l(k,n)^2}|D_{k,n}|\,\leq\,2\|\nu^n_0\|_{L^2}^2+2\sum_{k=1}^P\frac{(N^{k,n}_T)^2}{l(k,n)^2}|D_{k,n}|\,.
\end{eqnarray*}
Here the right hand side is independent of $t\leq T$ and thus we obtain
\begin{equation*}
\EX^n\sup\nolimits_{t\in[0,T]}\|\nu^n_t\|_{L^2}^2\,\leq\,2\EX^n\|\nu^n_0\|_{L^2}^2+2\sum_{k=1}^P\frac{\EX^n (N^{k,n}_T)^2}{l(k,n)^2}|D_{k,n}|\,\leq\,2\EX^n\|\nu^n_0\|_{L^2}^2+C_T\,,
\end{equation*}
where we have used that $N^{k,n}_T$ is Poisson distributed with rate $T\Lambda_{k,n}$ and thus $\EX^n(N^{k,n}_T)^2=T\Lambda_{k,n}+T^2\Lambda^2_{k,n}$. Here $C_T$ is some finite constant which depends on $T$ and the overall parameters of the model, i.e., $\tau,f,D$, but is independent of $k$ and $n$. Using this upper bound the triangle inequality yields the estimate
\begin{equation*}
\EX^n\sup\nolimits_{t\in[0,T]}\|\nu^n_t-\nu(t)\|_{L^2}^2\,\leq\, 2C_T\Bigl(\EX^n\|\nu^n_0\|_{L^2}^2+\|\nu\|_{C([0,T],L^2)}^2+1\Bigr)\,.
\end{equation*}
Therefore using the assumption $\sup_{n\in\mathbb{N}}\EX^n\|\nu_0^n\|_{L^2}^{2}<\infty$ it holds that
\begin{equation*}
\sup_{n\in\mathbb{N}}\EX^n\sup\nolimits_{t\in[0,T]}\|\nu^n_t-\nu(t)\|_{L^2}^{2} <\infty\,.
\end{equation*}

The general case for $r>1$ works analogously. Note that the $r$-th moment of the Poisson distribution is proportional to the $r$-th power of its rate. Hence, just as in the case of $r=1$, the term
\begin{equation*}
\EX^n\biggl(\sum_{k=1}^P\frac{(N^{k,n}_T)^{2}}{l(k,n)^2}|D_{k,n}|\biggr)^r
\end{equation*}
can thus be bounded from above by some constant $C_T$ independent of $k$ and $n$. The proof of Theorem \ref{lln} is completed.

\subsection{Proof of Corollary \ref{corollary_to_lln} (Corollary to the law of large numbers)}\label{section_proof_corol}

For $\alpha=0$ the statement of the corollary coincides with the statement of Theorem \ref{lln}, hence we consider $\alpha>0$. As in the proof of Theorem \ref{lln} we apply \cite[Thm.~4.1]{RTW} to the PDMPs $(Y^n_t)_{t\geq 0}$ however this time for the functions $\nu^n$ understood as taking values in the Hilbert space $H^{-\alpha}$ instead of $L^2$. Thus we have to validate again conditions (LLN1)--(LLN3) wherein the norm in $L^2$ is always replaced by the norm in $H^{-\alpha}$. The essential argument is sharpening the estimates in part (a) of the proof of Theorem \ref{lln} using optimal Sobolev Embedding Theorems such that the conditions \eqref{lln_sharpened_conditions} imply (LLN1). This we present in part (a) of the proof below. The Lipschitz condition (LLN2) of the Nemytzkii operator $F$ in the spaces $H^{-\alpha}$ is established in part (b). Finally, as the condition $\delta_+(n)\to\infty$ remains as in Theorem \ref{lln}, the condition (LLN3) follows immediately from the proof of Theorem \ref{lln} due to the continuous embedding of $L^2$ into $H^{-\alpha}$.\medskip

(a)\quad In the case $\alpha=0$, i.e., $H^\alpha=L^2$, we used in \eqref{lln_martingale_estimate} that $\|\mathbb{I}_{D_{k,n}}\|_{L^2}^2=|D_{k,n}|$. For general $\alpha>0$ we use the representation
\begin{equation*}
 \|\mathbb{I}_{D_{k,n}}\|_{H^{-\alpha}}=\sup_{\|\phi\|_{H^\alpha}}\bigl|(\phi,\mathbb{I}_{D_{k,n}})_{L^2}\bigr|.
\end{equation*}
In order to estimate the terms inside the supremum in the right hand side we use H\"older's inequality and the Sobolev embedding theorem, i.e., $H^\alpha(D)\hookrightarrow L^\infty(D)$ for $\alpha>d/2$ and $H^\alpha(D)\hookrightarrow L^r(D)$ with $r=d/(d/2-\alpha)$ for $0<\alpha<d/2$, see \cite[Thm.~7.34, Corol.~7.17]{Adams}. Thus we obtain
\begin{equation*}
\|\mathbb{I}_{D_{k,n}}\|_{H^{-\alpha}} \leq \left\{\begin{array}{cl}K_{d/(d/2-\alpha)}\,\|\mathbb{I}_{D_{k,n}}\|_{L^{\frac{2d}{d+2\alpha}}} & \textnormal{if } 0<\alpha<d/2, \\[3ex]
K_{\infty}\,\|\mathbb{I}_{D_{k,n}}\|_{L^1} & \textnormal{if } d/2<\alpha,
\end{array}\right.
\end{equation*}
where the constants $K$ are the constants arising from the continuous embeddings of the Sobolev spaces into the Lebesgue spaces. Evaluating the norms in the right hand side and further estimating using the maximal Lebesgue measure of the elements of the partition yields
\begin{equation*}
\|\mathbb{I}_{D_{k,n}}\|_{H^{-\alpha}}^2 \leq \left\{\begin{array}{cl}K_{d/(d/2-\alpha)}^2\,|D_{k,n}|\,v_+(n)^{2\alpha/d} & \textnormal{if } 0<\alpha<d/2, \\[3ex]
K_{\infty}^2\,|D_{k,n}|\,v_+(n) & \textnormal{if } d/2<\alpha,
\end{array}\right.
\end{equation*}

Note that the upper bounds are consistent with the condition in Theorem \ref{lln} for $\alpha=0$. Finally, as $H^{d/2}\hookrightarrow H^{(d/2-\eps)}$ for all small $\eps$, the result for $\alpha=d/2$ follows from the result above as
\begin{equation*}
\|\mathbb{I}_{D_{k,n}}\|_{H^{-d/2}}\,\leq\, \sup_{\|\phi\|_{H^{2\alpha}}}\bigl(\|\phi\|_{L^{d/\eps}}\,\|\mathbb{I}_{D_{k,n}}\|_{L^\frac{d}{d-\eps}}\bigr)\,\leq\,
C\,\|\mathbb{I}_{D_{k,n}}\|_{L^\frac{d}{d-\eps}}
\end{equation*}
where $C$ is the constant resulting from the continuous embedding of $H^{d/2}(D)$ into $H^{d/2-\eps}$. Thus we obtain for all $\eps>0$ the estimate
\begin{equation*}
\|\mathbb{I}_{D_{k,n}}\|^2_{H^{-d/2}}\,\leq\,C^2\,|D_{k,n}|\,v_+(n)^{\frac{d-2\eps}{d}}\,.  
\end{equation*}

(b)\quad Next we have to establish that the Nemytzkii operator $F$ on $L^2(D)$ is also Lipschitz continuous with respect to the norms $\|\cdot\|_{H^{-\alpha}}$, $\alpha\geq 0$, i.e., for all $\alpha\geq 0$ there exists a constant $L_{-\alpha}$ such that
\begin{equation}\label{lipschitz_condition_in_Hs}
\|F(g_1,t)-F(g_2,t)\|_{H^{-\alpha}}\,\leq \,L_{-\alpha}\,\|g_1-g_2\|_{H^{-\alpha}}\qquad \forall\, t\geq 0,\, g_1,g_2\in L^2(D)\,.
\end{equation}
We obtain due to the Lipschitz continuity of $f$, which implies absolute continuity of $f$, that
\begin{equation*}
\Big|\int_D \phi(x)\Bigl(F(g_1,t)(x)-F(g_2,t)(x)\Bigr)\,\dif x\Big| = \Big|\int_D \phi(x)\int_{z_1(t,x)}^{z_2(t,x)}f'(z)\,\dif z\,\dif x\Big|\,, 
\end{equation*}
where
\begin{equation*}
z_1(t,x)= \int_D w(x,y)g_1(y)\,\dif y+I(t,x), \qquad z_2(t,x)= \int_D w(x,y)g_2(y)\,\dif y+I(t,x)\,.
\end{equation*}
Applying H\"older's inequality and the essential boundedness of the derivative $f'$ we obtain the estimate
\begin{eqnarray*}
\Big|\int_D \phi(x)\Bigl(F(g_1,t)(x)-F(g_2,t)(x)\Bigr)\,\dif x\Big|&\leq& \|\phi\|_{L^p}\,\biggl(\int_D\biggl|\int_{z_1(t,x)}^{z_2(t,x)}f'(z)\,\dif z\biggr|^q\,\dif x\biggr)^{1/q}\\
&\leq& \|\phi\|_{L^p}\,\biggl(\int_D\biggl|\|f'\|_{L^\infty}\,\bigl(z_1(t,x)-z_2(t,x)\bigr)\biggr|^q\,\dif x\biggr)^{1/q}\\
&=& \|\phi\|_{L^p}\,\|f'\|_{L^\infty}\,\biggl(\int_D\biggl|\int_D w(x,y)\bigl(g_1(y)-g_1(y)\bigr)\,\dif y\biggr|^q\,\dif x\biggr)^{1/q}\,.
\end{eqnarray*}
Next, as by assumption $w(x,\cdot)\in H^{\alpha}$ we obtain
\begin{eqnarray*}
\biggl(\int_D\biggl|\int_D w(x,y)\bigl(g_1(y)-g_1(y)\bigr)\,\dif y\biggr|^q\,\dif x\biggr)^{1/q} &=& \biggl(\int_D\|w(x,\cdot)\|_{H^\alpha}^q\Big|\bigl\langle w(x,\cdot)/\|w(x,\cdot)\|_{H^\alpha}, g_1-g_2\bigr\rangle_{H^\alpha}\Big|^q\biggr)^{1/q}\\[2ex]
&\leq& \|w\|_{L^q\times H^\alpha}\, \|g_1-g_2\|_{H^{-\alpha}}\,.
\end{eqnarray*}
Overall this yields the estimate
\begin{equation*}
\Big|\bigl\langle\phi,F(g_1,t)-F(g_2,t)\bigr\rangle_{H^{-\alpha}}\Big|\,\leq\, \|\phi\|_{L^p}\,\|f'\|_{L^\infty}\,\|w\|_{L^q\times H^\alpha}\, \|g_1-g_2\|_{H^{-\alpha}}\,.
\end{equation*}
Hence taking the supremum on both sides of this inequality over all $\|\phi\|_{H^\alpha}=1$ we obtain the Lipschitz condition \eqref{lipschitz_condition_in_Hs} with $L_{-\alpha}:=L\,K_\alpha\,\|w\|_{L^q\times H^\alpha}$ where $K_\alpha$ is the constant resulting from the continuous embedding of $H^\alpha$ into $L^p$ and the Lipschitz constant $L$ of $f$ satisfies $L\geq \|f'\|_{L^\infty}$.

\subsection{Proof of Theorem \ref{theorem_infinite_time_lln} (Infinite time convergence)}\label{section_proof_infinite_time_lln}

(a)\, We first present an alternative representation for the jump processes $(\Theta^n_t)_{t\geq 0}$ and the solution $\nu$ of the Wilson-Cowan equation \eqref{Wilson_Cowan}. Using the generator of the PDMP $(\Theta^n_t,t)_{t\geq 0}$ we obtain that the components $\Theta^{k,n}$ satisfy
\begin{eqnarray}
 \Theta^{k,n}_t&=&\Theta^{k,n}_0+\int_0^t \lambda^n(\Theta^n_s,s)\int_{\mathbb{N}^p}\bigl(\xi^k-\Theta^{k,n}_s\bigr)\,\mu^n\bigl(\Theta^n_s,s;\dif\xi\bigr)\,\dif s +M^{k,n}_t\nonumber\\
&=& \Theta^{k,n}_0+\int_0^t\Bigl( -\frac{1}{\tau}\,\Theta^{k,n}_s+\frac{1}{\tau}\,l(k,n)\nbar f_{k,n}(\Theta^n_s,s)\Bigr)\,\dif s+M^{k,n}_t,\label{proof_2_SDE}
\end{eqnarray}
where $(M^{k,n}_t)_{t\geq 0}$ is a square-integrable c\`adl\`ag martingale given by
\begin{equation}\label{proof_2_def_martingale}
 M^{k,n}_t:=\Theta^{k,n}_t-\Theta^{k,n}_0 - \int_0^t \lambda^n(\Theta^n_s,s)\int_{\mathbb{N}^p}\bigl(\xi^k-\Theta^{k,n}_s\bigr)\,\mu^n\bigl(\Theta^n_s,s;\dif\xi\bigr)\,\dif s\,.
\end{equation}
As the jump process is regular this martingale is almost surely of finite variation and it could also be written in terms of a stochastic integral with respect to the associated martingale measure of the PDMP \cite{Jacobsen}. 
Next, considering $\Theta^{k,n}$ the solution of the stochastic differential equation \eqref{proof_2_SDE} driven by the martingale $M^{k,n}$ -- it is clear that a solution exists as the stochastic integral equation \eqref{proof_2_SDE} is constructed from a solution --  it follows from the variation of constants formula that it satisfies
\begin{equation}\label{proof_2_mild_solution}
\Theta^{k,n}_t=\e^{-t/\tau}\Theta^{k,n}_0+\frac{1}{\tau}\,l(k,n)\int_0^t \e^{-(t-s)/\tau}\,\nbar f_{k,n}(\Theta^n_s,s)\,\dif s+\int_0^t\e^{-(t-s)/\tau}\dif M^{k,n}_s\,. 
\end{equation}
This formula can also be easily verified path-by-path by inserting \eqref{proof_2_mild_solution} into \eqref{proof_2_SDE} and using integration by parts. Note that here the stochastic integral with respect to the martingale is just a Riemann-Stieltjes integral as the martingale is of finite variation. 
%
%
For the sake of completeness we briefly sketch the arguments. Thus, inserting \eqref{proof_2_mild_solution} into \eqref{proof_2_SDE} yields
\begin{eqnarray*}
 \Theta^{k,n}_t&=& \underbrace{\Theta^{k,n}_0-\frac{1}{\tau} \int_0^t\e^{-s/\tau}\Theta^{k,n}_0\,\dif s}_{(\ast)}\\
&&\mbox{} -\frac{1}{\tau}\,l(k,n)\biggl(\underbrace{\frac{1}{\tau}\int_0^t \int_0^s \e^{-(s-r)/\tau}\,\nbar f_{k,n}(\Theta^n_r,r)\,\dif r\,\dif s -\int_0^t\nbar f_{k,n}(\Theta^n_s,s)\,\dif s}_{(\ast\ast)}\biggr)\\
&&\mbox{} \underbrace{-\frac{1}{\tau}\,\int_0^t \int_0^s\e^{-(s-r)/\tau}\dif M^{k,n}_r\,\dif s+M^{k,n}_t}_{(\ast\ast\ast)}\,.
\end{eqnarray*}
Considering the three terms marked $(\ast)$ -- $(\ast\ast\ast)$ separately, we show that this right hand side equals \eqref{proof_2_mild_solution}. For the first term $(\ast)$ simply evaluating the integral yields
\begin{equation*}
\Theta^{k,n}_0-\frac{1}{\tau} \int_{0^t}\e^{-s/\tau}\Theta^{k,n}_0\,\dif s\, =\, \Theta^{k_n}_0-\frac{1}{\tau}\bigl(\alpha^{-1}\e^{-t/\tau}-\tau\bigr)\Theta^{k,n}_0\,=\,\e^{-t/\tau}\Theta^{k,n}_0\,, 
\end{equation*}
which gives the first term in the right hand side of \eqref{proof_2_mild_solution}. Next we simplify the term $(\ast\ast)$ employing integration by parts  to the first term in $(\ast\ast)$ which yields
\begin{eqnarray*}
\lefteqn{\frac{1}{\tau}\int_0^t \int_0^s \e^{-(s-r)/\tau}\,\nbar f_{k,n}(\Theta^n_r,r)\,\dif r\,\dif s\ =}\\[2ex]
&&\phantom{xxxxxxxxxxxxxxxxx}=\ \frac{1}{\tau}\int_0^t \e^{(t-s)/\tau} \int_0^s \e^{-(t-r)/\tau}\,\nbar f_{k,n}(\Theta^n_r,r)\,\dif r\,\dif s\\
&&\phantom{xxxxxxxxxxxxxxxxx}=\ \frac{1}{\tau}\,\biggl(-\tau\,\e^{(t-s)/\tau}\int_0^s\e^{-(t-r)/\tau}\,\nbar f_{k,n}(\Theta^n_r,r)\,\dif r\biggr)\bigg|_0^t\\ 
&&\phantom{xxxxxxxxxxxxxxxxx=}\ \mbox{}-\frac{1}{\tau}\int_0^t\bigl(-\tau\bigr)\e^{(t-s)\tau}\e^{-(t-s)/\tau}\,\nbar f_{k,n}(\Theta^n_s,s)\,\dif s\\
&&\phantom{xxxxxxxxxxxxxxxxx}=\ -\int_0^t\e^{-(t-s)/\tau}\,\nbar f_{k,n}(\Theta^n_s,s)\,\dif s +\int_0^t\nbar f_{k,n}(\Theta^n_s,s)\,\dif s\,.
\end{eqnarray*}
Thus we obtain subtracting from this right hand side the second term in $(\ast\ast)$ that
\begin{equation*}
(\ast\ast)\,=\, -\int_0^t\e^{-(t-s)/\tau}\,\nbar f_{k,n}(\Theta^n_s,s)\,\dif s\,.
\end{equation*}
This term is just the second term in the right hand side of \eqref{proof_2_mild_solution}. It remains to consider the term marked $(\ast\ast\ast)$. We have already stated that the stochastic integral with respect to the martingale \eqref{proof_2_def_martingale} is defined path-by-path as a Riemann-Stieltjes integral and thus satisfies
\begin{eqnarray}
\lefteqn{\hspace{-1.0ex}-\frac{1}{\tau}\int_0^s\e^{-(s-r)/\tau}\dif M^{k,n}_s\ =}\label{proof_2_Riemann_Stieltjes}\\[2ex]
&=&-\frac{1}{\tau}\sum_{\tau^n_j\leq s} \e^{-(s-\tau^n_j)/\tau}\bigl(\Theta^{k,n}_{\tau^n_j}-\Theta^{k,n}_{\tau^n_{j-}}\bigr)+\frac{1}{\tau}\int_0^s\e^{-(s-r)/\tau}\lambda^n(\Theta^n_r,r)\int_{\mathbb{N}_0^P}\bigl(\xi^k-\Theta^{k,n}_r\bigr)\,\mu^n\bigl((\Theta^n_r,r),\dif\xi)\,\dif s\,,\nonumber
\end{eqnarray}
where $\tau^n_j$ denotes the $j$-th jump time of the $n$-th PDMP. Integrating the sum in this right hand side over $(0,t)$ yields
\begin{eqnarray*}
 -\frac{1}{\tau}\int_0^t\sum_{\tau^n_j\leq s} \e^{-(s-\tau^n_j)/\tau}\bigl(\Theta^{k,n}_{\tau^n_j}-\Theta^{k,n}_{\tau^n_{j-}}\bigr)\,\dif s &=& \sum_{\tau_j^n\leq t}\e^{-(t-\tau_j^n)/\tau}\bigl(\Theta^{k,n}_{\tau^n_j}-\Theta^{k,n}_{\tau^n_{j-}}\bigr) - \sum_{\tau_j^n\leq t}\bigl(\Theta^{k,n}_{\tau^n_j}-\Theta^{k,n}_{\tau^n_{j-}}\bigr)\\
&=& \sum_{\tau_j^n\leq t}\e^{-(t-\tau_j^n)/\tau}\bigl(\Theta^{k,n}_{\tau^n_j}-\Theta^{k,n}_{\tau^n_{j-}}\bigr) -\bigl(\Theta^{k,n}_t-\Theta^{k,n}_0\bigr)\,.
\end{eqnarray*}
Next, we apply integration by parts to the integral over $(0,t)$ of the second term above analogously to the application to term $(\ast\ast)$ and obtain
\begin{eqnarray*}
\lefteqn{\frac{1}{\tau}\int_0^t \int_0^s\e^{-\alpha(s-r)}\lambda^n(\Theta^n_r,r)\int_{\mathbb{N}_0^P}\bigl(\xi^k-\Theta^{k,n}_r\bigr)\,\mu^n\bigl((\Theta^n_r,r),\dif\xi)\,\dif r\,\dif s}\\[2ex]
&=& -\int_0^t\e^{-(t-s)/\tau}\lambda^n(\Theta^n_s,s)\int_{\mathbb{N}_0^P}\bigl(\xi^k-\Theta^{k,n}_s\bigr)\,\mu^n\bigl((\Theta^n_s,s),\dif\xi)\,\dif s\\
&&\phantom{xxxxxxxxxxxxxxxxxxxxxxxxxxxxxxxxxxxxx}\mbox{} +\int_0^t\lambda^n(\Theta^n_s,s)\int_{\mathbb{N}_0^P}\bigl(\xi^k-\Theta^{k,n}_s\bigr)\,\mu^n\bigl((\Theta^n_s,s),\dif\xi)\,\dif s\,.
\end{eqnarray*}

Hence, overall these considerations show that
\begin{equation*}
(\ast\ast\ast)\,=\, \int_0^t\e^{-(t-s)/\tau}\,\dif M^{k,n}_s
\end{equation*}
and we obtain the final, third term in the right hand side of \eqref{proof_2_mild_solution}. This completes the proof that \eqref{proof_2_mild_solution} solves the equation \eqref{proof_2_SDE}. 
\medskip

Further, we obtain from the variation of constants formula for $\Theta^{k,n}_t$ also a representation for the stochastic mean activity $\nu^n$ by inserting \eqref{proof_2_mild_solution} into its definition \eqref{def_coordinate_function}. This gives
\begin{eqnarray}
\nu^n_t&=&\e^{-t/\tau}\nu^n_0+\frac{1}{\tau}\,\sum_{k=1}^P\int_0^t \e^{-(t-s)/\tau}\,\nbar f_{k,n}(\Theta^n_s,s)\,\dif s\,\mathbb{I}_{D_{k,n}} +\sum_{k=1}^P\frac{1}{l(k,n)}\int_0^t\e^{-(t-s)/\tau}\dif M^{k,n}_s\, \mathbb{I}_{D_{k,n}} \nonumber\\    
&=&\e^{-t/\tau}\nu^n_0+\frac{1}{\tau}\int_0^t \e^{-(t-s)/\tau}\,\nbar F^n(\nu^n_s,s)\,\dif s +\sum_{k=1}^P\frac{1}{l(k,n)}\int_0^t\e^{-(t-s)/\tau}\dif M^{k,n}_s\, \mathbb{I}_{D_{k,n}}\,.\label{stoch_sol_var_of_constants}
\end{eqnarray}

Finally, in order to compare stochastic and deterministic solutions we use that the solution of the Wilson-Cowan equation can also be given via the variation of constants formula, i.e., it holds that for all $t\geq 0$
\begin{equation}\label{det_sol_var_of_constants}
 \nu(t)=\e^{-t/\tau}\nu(0)+\frac{1}{\tau}\int_0^t\e^{-(t-s)/\tau}\,F(\nu(s),s)\,\dif s\,.
\end{equation}

Thus, subtracting \eqref{det_sol_var_of_constants} from \eqref{stoch_sol_var_of_constants} and taking the expectation of the norm in $H^{-\alpha}$ yields the estimate
\begin{eqnarray}
\EX^n\|\nu(t)-\nu^n_t\|_{H^{-\alpha}}&=&\e^{-t/\tau}\,\EX^n\|\nu(0)-\nu_0^n\|_{H^{-\alpha}}+\frac{1}{\tau}\int_0^t\e^{-(t-s)/\tau}\EX^n\bigr\|F(\nu(s),s)-\nbar F^n(\nu^n_s,s)\bigl\|_{H^{-\alpha}}\nonumber\\
&&\mbox{}+\EX^n\Big\|\sum_{k=1}^P\frac{1}{l(k,n)}\int_0^t\e^{-(t-s)/\tau}\dif M^{k,n}_s\, \mathbb{I}_{D_{k,n}}\Big\|_{H^{-\alpha}}\,.\label{infinite_time_convergence_gronwall_est_1}
\end{eqnarray}
We deal with the terms in the right hand side of \eqref{infinite_time_convergence_gronwall_est_1} separately in the following such that we can apply Gronwall's inequality. Note that the term containing the initial condition vanishes due to the assumptions of the theorem. We start with the stochastic integrals in the subsequent part (b) of the proof.

\medskip

(b)\, As due to Jensen's inequality $\EX|Y|\leq\sqrt{\EX|Y|^2}$ it makes sense to calculate the second moment of the stochastic integral in the right hand side. For the norm in $H^{-\alpha}$ we use $\|\phi\|_{H^{-\alpha}}^2=(\phi,\phi)_{H^{-\alpha}}$ and thus obtain using the linearity of the inner product
\begin{eqnarray*}
\Big\|\sum_{k=1}^P\frac{1}{l(k,n)}\underbrace{\int_0^t\e^{-(t-s)/\tau}\dif M^{k,n}_s}_{=:\beta^{k,n}_t}\, \mathbb{I}_{D_{k,n}}\Big\|_{H^{-\alpha}}^2
&=&\sum_{k=1}^P\frac{|\beta_{k,n}|^2}{l(k,n)^2}\,\|\mathbb{I}_{D_{k,n}}\|_{H^{-\alpha}}^2\\[1ex]
&&\mbox{}+ \sum_{\substack{k,j=1\\ k\neq j}}^P\frac{\beta_{k,n}\,\beta_{j,n}}{l(k,n)\,l(j,n)}\,(\mathbb{I}_{D_{k,n}},\mathbb{I}_{D_{j,n}})_{H^{-\alpha}}\,.
\end{eqnarray*}
We next consider the individual expectations of the random terms $|\beta_{k,n}|^2$ and $\beta_{k,n}\,\beta_{j,n}$ in the right hand side. 
We have already stated that the stochastic integral with respect to the martingale \eqref{proof_2_def_martingale} is defined path-by-path as a Riemann-Stieltjes integral, see \eqref{proof_2_Riemann_Stieltjes} and, moreover, \eqref{proof_2_Riemann_Stieltjes} implies that the stochastic convolution integral can be written as a stochastic integral with respect to the fundamental martingale measure $M^n$ associated with the PDMP $(\Theta^n_t,t)_{t\geq 0}$, see \cite{Jacobsen}, i.e.,
\begin{eqnarray*}
\int_0^t\e^{-(t-s)/\tau}\dif M^{k,n}_s&=& \int_{[0,t]\times {\mathbb{N}^P_0}} \e^{-(t-s)/\tau}\bigl(\xi^k-\Theta^{k,n}_{s-}\bigr)\,M^n(\dif s,\dif\xi)
\end{eqnarray*}
with predictable integrand
\begin{equation*}
(\xi,s,\omega)\mapsto \e^{-(t-s)/\tau}\bigl(\xi^k-\Theta^{k,n}_{s-}(\omega)\bigr)\,.
\end{equation*}
Then we obtain due to the It\^o-isometry following from \cite[Prop.~4.6.2]{Jacobsen} using \eqref{components_mean_bound} that
\begin{eqnarray*}
\EX^n\Big|\int_0^t\e^{-(t-s)/\tau}\dif M^{k,n}_s\Big|^2 &=& \EX^n \int_0^t\lambda^n(Y^n_s)\int_{\mathbb{N}^P_0}\e^{-2(t-s)/\tau}\bigl(\xi^k-\Theta^{k,n}_s\bigr)^2\mu^n(Y^n_s,\dif\xi)\,\dif s \\
&\leq& \EX^n \int_0^t\e^{-2(t-s)/\tau}\,\Bigl(\frac{1}{\tau}\,\Theta^{k,n}_s+\frac{1}{\tau}\,l(k,n)\nbar f_{k,n}(Y^n_s)\Bigr)\,\dif s \\
&\leq& \frac{1}{\tau}\Bigl(\EX^n\Theta^{k,n}_0+2\,l(k,n)\,\|f\|_0\Bigr)\int_0^t\e^{-2(t-s)/\tau}\dif s\,.
\end{eqnarray*}
It remains to consider the product $\beta^{k,n}_t\,\beta^{f,n}_t$ for which we obtain due to the integration by parts formula
\begin{equation}
\beta^{k,n}_t\,\beta^{j,n}_t\,=\,\int_0^t \beta^{k,n}_{s-}\,\dif \beta^{j,n}_s + \int_0^t \beta^{j,n}_{s-}\,\dif \beta^{k,n}_s+\bigl[\beta^{k,n},\beta^{j,n}\bigr]_t\,,
\end{equation}
where the square brackets denote the quadratic variation process. The expectation of each of the terms in the right hand side vanishes: The first two are stochastic integrals with respect to martingales, hence martingales themselves which are identical to zero at the origin. Furthermore, as both martingales are c\`adl\`ag with paths of finite variation on compacts, hence quadratic pure jump martingales, we obtain for the quadratic variation process
\begin{eqnarray*}
\bigl[\beta^{k,n},\beta^{j,n}\bigr]_t&=&\sum_{s\leq t}\bigl(\beta^{k,n}_s-\beta^{k,n}_{s-}\bigr)\bigl(\beta^{j,n}_s-\beta^{j,n}_{s-}\bigr)\,.
\end{eqnarray*}
However, as all jump times of the two martingales a.s.~differ it follows that $\bigl[\beta^{k,n},\beta^{j,n}\bigr]_t=0$.

\medskip

Thus overall we have established that
\begin{eqnarray}
\EX^n\Big\|\sum_{k=1}^P\frac{1}{l(k,n)}\int_0^t\e^{-(t-s)/\tau}\dif M^{k,n}_s\, \mathbb{I}_{D_{k,n}}\Big\|_{H^{-\alpha}}^2&=&
\frac{1}{2}\,\sum_{k=1}^P\frac{1+2\,\|f\|_0}{l(k,n)}\,\|\mathbb{I}_{D_{k,n}}\|_{H^{-\alpha}}^2,
\end{eqnarray}
where $1/2$ is an upper bound for $\frac{1}{\tau}\int_0^t\e^{-2(t-s)/\tau}\dif s$ independent of $t$. Estimating the norm $\|\mathbb{I}_{D_{k,n}}\|_{H^{-\alpha}}^2$ just as in the proof of Corollary \ref{corollary_to_lln} we finally obtain that
\begin{equation*}
\EX^n\Big\|\sum_{k=1}^P\frac{1}{l(k,n)}\int_0^t\e^{-(t-s)/\tau}\dif M^{k,n}_s\, \mathbb{I}_{D_{k,n}}\Big\|_{H^{-\alpha}}\,\leq\,\frac{1}{4}\biggl(\bigl(1+2\|f\|_0\bigr)\,|D|\,\frac{v_+(n)^r}{\ell_-(n)}\biggr)^{1/2},
\end{equation*}
with $r=2\alpha/d$ for $0\leq \alpha\leq d/2$, $r=1-\eps$ for $\alpha=d/2$ and $r=1$ for $\alpha>d/2$.

\medskip

(c)\quad We next estimate the term
\begin{equation*}
\int_0^t\e^{-(t-s)/\tau}\EX^n\bigr\|F(\nu(s),s)-\nbar F^n(\nu^n_s,s)\bigl\|_{H^{-\alpha}}
\end{equation*}
in \eqref{infinite_time_convergence_gronwall_est_1}. From part (b) of the proof of Theorem \ref{lln} in Section \ref{section_proof_lln} it follows that
\begin{equation}\label{infin_lln_gan_conv}
\EX^n\|F(\nu^n_t,t)-\nbar F^n(\nu^n_t,t)\|_{H^{-\alpha}}\,\leq\, \delta_+(n)\,\frac{K_{-\alpha}\,L}{\pi}\Bigl(\sqrt{|D|}\,(1+\|f\|_0)\, \|\nabla_{\!x} w\|_{L^2\times L^2}+\|\nabla_{\!x} I(t)\|_{L^2}\Bigr)\,,
\end{equation}
where $F$ is the Nemyztkii operator defined in \eqref{def_Nemytzkii_op} and $K_{-\alpha}$ is a constant resulting from the continuous embedding of $L^2$ into $H^{-\alpha}$. Here, the right hand side can be further estimated independently of $t\geq 0$ using the assumption that $\|\nabla_{\!x} I(t)\|_{L^2}$ is uniformly bounded in $t\geq 0$. Furthermore we have shown in Section \ref{section_proof_corol} in the proof of Corollary \ref{corollary_to_lln}, that under the appropriate assumptions the Nemytzkii operator $F$ is Lipschitz continuous on $H^{-\alpha}$, $\alpha\geq 0$, with Lipschitz constant $L_{-\alpha}>0$ independent of $t\geq0$, i.e.,
\begin{equation}\label{infin_lln_lipschitz_2}
\|F(g_1,t)-F(g_2,t)\|_{H^{-\alpha}}\,\leq\, L_{-\alpha}\,\|g_1-g_2\|_{H^{-\alpha}}\qquad\forall\,g_1,g_2\in L^2\,.
\end{equation}
A combination of the triangle inequality and the estimates \eqref{infin_lln_gan_conv} and \eqref{infin_lln_lipschitz_2} yields
\begin{eqnarray*}
\int_0^t\e^{-(t-s)/\tau}\EX^n\bigr\|F(\nu(s),s)-\nbar F^n(\nu^n_s,s)\bigl\|_{H^{-\alpha}}&\leq& L_{-\alpha}\int_0^t\e^{-(t-s/\tau)}\EX^n\|\nu(s)-\nu^n_s\|_{H^{-\alpha}}\dif s +\landau\bigl(\delta_+(n)\bigr)\,.
\end{eqnarray*}\medskip

Overall, it thus follows from \eqref{infinite_time_convergence_gronwall_est_1} that
\begin{eqnarray*}
\EX^n\|\nu(t)-\nu^n_t\|_{H^{-\alpha}}&\leq&\EX^n\|\nu(0)-\nu_0^n\|_{H^{-\alpha}}+\frac{L_{-\alpha}}{\tau}\int_0^t\e^{-(t-s/\tau)}\EX^n\|\nu(s)-\nu^n_s\|_{H^{-\alpha}}\dif s\\
&&\mbox +\landau\biggl(\delta_+(n)+\sqrt{\frac{v_+(n)^r}{\ell_-(n)}}\biggr)\,.
\end{eqnarray*}
Then an application of Gronwall's inequality yields
\begin{eqnarray*}
\EX^n\|\nu(t)-\nu^n_t\|_{H^{-\alpha}}&\leq& \biggl(\EX^n\|\nu(0)-\nu_0^n\|_{H^{-\alpha}}+\landau\biggl(\delta_+(n)+\sqrt{\frac{v_+(n)^r}{\ell_-(n)}}\biggr)\biggr)\exp\Bigl(\frac{L_{-\alpha}}{\tau}\int_0^t\e^{-(t-s)/\tau}\dif s\Bigr)\\
&\leq& \biggl(\EX^n\|\nu(0)-\nu_0^n\|_{H^{-\alpha}}+\landau\biggl(\delta_+(n)+\sqrt{\frac{v_+(n)^r}{\ell_-(n)}}\biggr)\biggr)\,\e^{L_{-\alpha}}\,.
\end{eqnarray*}
By assumptions of the theorem the term in the right hand side converges to zero for $n\to\infty$. As this convergence is uniform in $t$ it holds that
\begin{equation}
\lim_{n\to\infty} \sup_{t\geq 0}\EX^n\|\nu(t)-\nu^n_t\|_{H^{-\alpha}} =0\,.
\end{equation}

\subsection{Proof of Theorem \ref{martingale_clt} (Martingale central limit theorem)}\label{section_proof_martingale_clt}

In order to prove the martingale central limit theorem we employ the general martingale central limit theorem \cite[Thm.~5.1]{RTW} for the Hilbert space $H^{-\alpha}$, i.e., the dual of the Sobolev space $H^{\alpha}$, for $\alpha>d$. To apply this theorem it suffices to prove the following conditions. Subsequently we use $\rho_n=\sqrt{\ell_-(n)/v_+(n)}$ to denote the rescaling sequence and use the notation
\begin{eqnarray}\label{jump_quadratic_var}
\langle G^n(t)\phi,\phi\rangle_{H^{\alpha}}&=&\lambda(Y^n_t)\int_{\mathbb{N}^P_0} \langle\nu^n(\xi)-\nu^n(\Theta^n_t),\phi\rangle^2_{H^{\alpha}}\,\mu^n(Y^n_t,\dif\xi)
\end{eqnarray}
which corresponds to the quadratic variation of the martingales $(M^n_t)_{t\geq 0}$, see \cite{RTW} for a discussion.

\setlength{\leftmargini}{4.5em}
\begin{itemize}[align=left]
  \setlength{\labelwidth}{4.0em}
  \setlength{\labelsep}{0.5em}
\item[\textnormal{(CLT1)}] For all $t>0$ it holds that
\begin{equation}\label{first_clt_section_theorem_2nd_mom_cond}
\sup_{n\in\mathbb{N}}\rho_n\,\EX^n\int_0^t\Bigl[\lambda^n(Y^n_s)\int_{\mathbb{N}^P}\|\nu^n(\xi)-\nu^n(\Theta^n_s)\|_{H^{-\alpha}}^2\,\mu^n\bigl(Y^n_s,\dif\xi\bigr)\,\dif s\Bigr]\,<\infty\,,
\end{equation}
and there exists an orthonormal basis $(\varphi_j)_{j\in\mathbb{N}}$ of $H^{\alpha}(D)$ such that for all $j\in\mathbb{N}$ and $t\geq 0$
\begin{equation}\label{alternative_tightness_cond}
\rho_n\,\EX^n\langle G^n(Y^n_t)\varphi_j, \varphi_j\rangle_{H^{\alpha}}\ \leq\  \gamma_j \, C\,,
\end{equation}
where the constants $\gamma_j>0$ are independent of $n$ and $t$, satisfy $\sum_{j\geq 1}\gamma_j<\infty$, and the constant $C>0$ is independent of $n$ and $k$ but may depend on $t$.

\item[\textnormal{(CLT2)}] The jump heights of the rescaled martingales are almost surely uniformly bounded, i.e., there exists a constant $\beta<\infty$ such that it holds almost surely for all $n\in\mathbb{N}$ that
\begin{equation}\label{clt_martingale_cond_1}
\sup_{t\geq 0}\,\sqrt{\rho_n}\,\big\|\nu^n(\Theta^n_t)-\nu^n(\Theta^n_{t-})\big\|_{H^{-\alpha}}<\beta\,.
\end{equation}

Further, for all $\phi\in H^{\alpha}$ and all $t>0$ it holds that
\begin{equation}\label{clt_martingale_cond_2}
\lim_{n\to\infty} \int_0^t\EX^n\big|\bigl\langle G(\nu(s))\,\phi,\phi\bigr\rangle_{H^{\alpha}}-\rho_n\bigl\langle G^n(Y^n_s)\,\phi,\phi\bigr\rangle_{H^{\alpha}}\big|\,\dif s\, =\, 0\,.
\end{equation}

\end{itemize}

On a technical level we note that the condition (CLT1) guarantees tightness of the sequence of rescaled martingales $(\sqrt{\rho_n}\,M^n_t)_{t\geq 0}$ in the Skorokhod space of c\`adl\`ag functions in $H^{-\alpha}$. This property is equivalent to relative compactness in the topology of weak convergence of measures and thus implies the existence of a convergent subsequence. The conditions (CLT2) are then sufficient to establish that any limit possesses the form of a diffusion process defined by the covariance operator $C$ given in \eqref{def_covariance_op}. In particular, condition \eqref{clt_martingale_cond_2} precisely gives the convergence of the quadratic variations and is thus the central condition. In the subsequent two parts of the proof we show that they are satisfied: In part (a) we prove conditions \eqref{first_clt_section_theorem_2nd_mom_cond} and \eqref{alternative_tightness_cond}  and part (b) establishes \eqref{clt_martingale_cond_1} and \eqref{clt_martingale_cond_2}.\medskip

(a)\quad We first prove conditions \eqref{first_clt_section_theorem_2nd_mom_cond} and \eqref{alternative_tightness_cond}. Here we also observe the significance of the choice of the norm in $H^{-\alpha}$ with $\alpha>d$ for establishing the convergence, which is essentially that it guarantees the existence a Sobolev space $H^{\alpha_1}$ with continuous embeddings $H^{\alpha}\hookrightarrow H^{\alpha_1}\hookrightarrow C(\nbar D)$, where the first is of Hilbert-Schmidt type. 
For subsequent use we recall the estimates 
\begin{equation*}
\|\mathbb{I}_{D_{k,n}}\|_{H^{-\alpha}}^2\leq K_\alpha^2\,|D_{k,n}|^2
\end{equation*}
with a suitable constant $K_\alpha>0$, which we have already established in the proof of Corollary \ref{corollary_to_lln} due to the H\"older inequality and the Sobolev Embedding Theorem. 
Therefore we obtain for the term inside the expectation in \eqref{first_clt_section_theorem_2nd_mom_cond} the estimate
\begin{equation*}
\lambda^n(Y^n_s)\int_{\mathbb{N}^P_0}\|\nu^n(\xi)-\nu^n(\Theta^n_s)\|_{H^{-\alpha}}^2\,\mu^n\bigl(Y^n_s,\dif\xi\bigr)\,\leq\, \frac{1}{\tau}\,K_{\alpha}^2\sum_{k=1}^P \frac{|D_{k,n}|^{2}}{l(k,n)^2}\Bigl(\Theta^{k,n}_s+l(k,n)\,\nbar f_{k,n}(Y^n_s)\Bigr)\,.
\end{equation*}
Next taking the expectation, using the bound \eqref{components_mean_bound} on $\EX^n\Theta^{k,n}_s$ and integrating over $[0,t]$ we obtain the estimate
\begin{equation*}
\int_0^t\EX^n\Bigl[\lambda^n(Y^n_s)\int_{\mathbb{N}^P_0}\|\nu^n(\xi)-\nu^n(\Theta^n_s)\|_{H^{-\alpha}}^2\,\mu^n\bigl(Y^n_s,\dif\xi\bigr)\Bigr]\,\dif s\,\leq\, \frac{t}{\tau}\,K_{\alpha}^2(1+2\|f\|_0)\, \frac{v_+(n)}{\ell_-(n)}\,.
\end{equation*}
Multiplying both sides with $\rho_n=\ell_-(n)/v_+(n)$ we find that condition \eqref{first_clt_section_theorem_2nd_mom_cond} is satisfied.\medskip

We proceed to condition \eqref{alternative_tightness_cond} and first of all expand the integrand to obtain
\begin{eqnarray*}
\langle G^n(Y^n_s)\varphi_j, \varphi_j\rangle_{H^{\alpha}}&=& \lambda^n(Y^n_s)\int_{\mathbb{N}^P_0}\bigl\langle \nu^n(\xi)-\nu^n(\Theta^n_s),\phi_j\bigr\rangle_{H^{-\alpha}}^2\,\mu^n\bigl(Y^n_s,\dif\xi\bigr)\\
&=& \frac{1}{\tau}\sum_{k=1}^P \frac{1}{l(k,n)^2}\Bigl(\Theta^{k,n}_s+l(k,n)\,\nbar f_{k,n}(Y^n_s)\Bigr)\,\langle\mathbb{I}_{D_{k,n}},\varphi_j\rangle_{H^{\alpha}}^2\,.
\end{eqnarray*}
We next estimate the term $\langle\mathbb{I}_{D_{k,n}},\varphi\rangle_{H^{\alpha}}^2$. Here we use the fact that for a function in $L^2(D)$ its application as an element of the dual $H^{-\alpha}$ as well as $H^{-\alpha_1}$ for any $\alpha_1$ with $0\leq \alpha_1<\alpha$ coincide. We choose $\alpha_1$ such that $d/2<\alpha_1<\alpha-d/2$ and obtain
\begin{equation*}
\langle\mathbb{I}_{D_{k,n}},\varphi\rangle_{H^{\alpha}}^2\leq \|\mathbb{I}_{D_{k,n}}\|_{H^{-\alpha_1}}^2\|\varphi_j\|_{H^{\alpha_1}}^2\leq K_{\alpha_1}^2\,|D_{k,n}|^2\,\|\varphi_j\|_{H^{\alpha_1}}^2,
\end{equation*}
where $K_{\alpha_1}$ is the constant resulting from the Sobolev Embedding Theorem. Next taking the expectation, estimating the expectation terms as before and multiplying by $\rho_n$ yields
\begin{equation*}
\rho_n\,\EX^n\langle G^n(Y^n_t)\varphi_j, \varphi_j\rangle_{H^{\alpha}}\ \leq\ \frac{1}{\tau}\,K_{\alpha_1}^2(1+\|f\|_0)\,\|\varphi_j\|_{H^{\alpha_1}}^2\,.
\end{equation*}
We chose the constants in \eqref{alternative_tightness_cond} as $C:=K_{\alpha_1}^2(1+\|f\|_0)/\tau$ and $\gamma_j:=\|\varphi_j\|_{H^{\alpha_1}}^2$. Finally, as due to Maurin's Theorem the embedding of the space $H^\alpha$ into $H^{\alpha_1}$ is of Hilbert-Schmidt type, cf.~footnote \ref{maurin_footnote} on p.~\pageref{maurin_footnote}, it holds that $\sum_{j\geq 1}\|\varphi_j\|^2_{H^{\alpha_1}}<\infty$. Condition \eqref{alternative_tightness_cond} is satisfied.
\medskip

(b)\quad The estimates in part (a) further show that the jump sizes are almost surely uniformly bounded as
\begin{equation*}
\sup_{t\geq 0}\sqrt{\rho_n}\,\big\|\nu^n(\Theta^n_t)-\nu^n(\Theta^n_{t-}) \big\|_{H^{-\alpha}}\leq \,K_{\alpha}\sqrt{\frac{v_+(n)}{\ell_-(n)}}\,.
\end{equation*}
Here the upperbound in the right hand side converges to zero for $n\to\infty$ and thus the left hand side is bounded over all $n\in\mathbb{N}$. Therefore condition \eqref{clt_martingale_cond_1} holds and we are left to prove the convergence of the quadratic variation \eqref{clt_martingale_cond_2}. For the jump process the quadratic variation satisfies
\begin{eqnarray*}
\langle G^n(t)\phi,\phi\rangle_{H^{\alpha}}&=&\lambda(Y^n_t)\int_{\mathbb{N}^P_0} \langle\nu^n(\xi)-\nu^n(\Theta^n_t),\phi\rangle^2_{H^{\alpha}}\,\mu^n(Y^n_t,\dif\xi)\\ 
&=& \frac{1}{\tau}\sum_{k=1}^P \frac{1}{l(k,n)^2}\bigl(\Theta^{k,n}_t+l(k,n) \nbar f_{k,n}(Y^n_t)\,\langle\mathbb{I}_{D_{k,n}},\phi\rangle_{H^{\alpha}}^2\,.  
\end{eqnarray*}
The quadratic variation of the limiting diffusion is given by
\begin{equation*}
\langle G(\nu(t),t)\phi,\phi\rangle_{H^{\alpha}} \,=\, \int_D \phi(x)^2\Bigl(\frac{1}{\tau}\,\nu(t,x)+\frac{1}{\tau}\,f\Bigl(\int_Dw(x,y)\nu(t,y)\,\dif y+I(t,x)\Bigr)\,\dif x\,.
\end{equation*}

Here the necessary estimates are split into several parts which are separately considered in the following. Afterwards, the estimates are combined to infer the convergence \eqref{clt_martingale_cond_2}. In the following we use again $F$ as the Nemytzkii operator defined in \eqref{def_Nemytzkii_op}. 
Hence, for the difference of the quadratic variations we obtain the estimate
\begin{eqnarray}
\lefteqn{\hspace{-25pt}\EX^n\big|\langle G(\nu(t),t)\phi,\phi\rangle_{H^{\alpha}}-\rho_n\langle G^n(t)\phi,\phi\rangle_{H^\alpha}\big|}\nonumber\\[2ex]
&&=\ \frac{1}{\tau}\,\EX^n\Big|\int_D \phi(x)^2\nu(t,x)+\phi(x)^2F(\nu(t),t)(x)\,\dif x\nonumber\\
&&\phantom{xxxxxxxxxxxxxxxxxxxxxixxxxx}\mbox{}-\sum_{k=1}^P \frac{\rho_n}{l(k,n)^2}\bigl(\Theta^{k,n}_t+l(k,n)\,\nbar f_{k,n}(Y^n_t)\,\langle\mathbb{I}_{D_{k,n}},\phi\rangle_{H^{\alpha}}^2\Big|\nonumber\\[2ex]
&&\leq\ \frac{1}{\tau}\,\EX^n\Big|\int_D \underbrace{\phi(x)^2\nu(t,x)}_{(i)}+\underbrace{\phi(x)^2 F(\nu(t),t)(x)}_{(ii)}\,\dif x\nonumber\\
&&\phantom{xxxxxxxxxxxxxxxxxxxxxixx}\mbox{} - \int_D \underbrace{\phi(x)^2 \nu^n(\Theta^n_t)(x)}_{(i)}+\underbrace{\phi(x)^2F(\nu^n(\Theta^n_t),t)(x)}_{(ii)}\,\dif x\,\Big|\nonumber\\[2ex]
&&\phantom{\leq}\ \mbox{}+\frac{1}{\tau} \,\EX^n\Big|\int_D \underbrace{\phi(x)^2 \nu^n(\Theta^n_t)(x)}_{(iii)}+\underbrace{\phi(x)^2F(\nu^n(\Theta^n_t),t)(x)}_{(iv)}\,\dif x\label{the_basic_long_estimate_quad_var_proof}\\
&&\phantom{xxxxxxxxxixxxxxx}\mbox{}-\sum_{k=1}^P\frac{\rho_n}{l(k,n)} \Bigl(\underbrace{\frac{\Theta^{k,n}_t}{l(k,n)}\,\langle\mathbb{I}_{D_{k,n}},\phi\rangle_{H^{\alpha}}^2}_{(iii)}+\underbrace{\frac{l(k,n)}{l(k,n)}\nbar f_{k,n}(Y^n_t)\,\langle\mathbb{I}_{D_{k,n}},\phi\rangle_{H^\alpha}^2}_{(iv)}\Bigr)\Big|\,.\nonumber
\end{eqnarray}
Using the triangle inequality once again for each of the two differences grouping the terms marked ($i$)--($iv$) we obtain four terms which we subsequently estimate separately. Finally, in part ($v$) we combine the four estimates.\medskip

($i$)\quad The first term is the simplest to estimate. Using the Cauchy-Schwarz inequality we obtain 
\begin{equation}\label{quad_var_proof_est1}
\EX^n\Big|\int_D \phi^2(x)\bigl(\nu(t,x)-\nu^n_t(x)\bigr)\,\dif x\Big| \,\leq\, \|\phi\|_{L^4}^2\,\EX^n\|\nu(t)-\nu^n_t\|_{L^2}\,.
\end{equation}

\medskip

($ii$)\quad We next consider the difference arising from the terms marked ($ii$) and obtain using the Lipschitz condition \eqref{lipschitz_on_f} on $f$ and the Cauchy-Schwarz inequality twice
\begin{eqnarray}
\EX^n\Big|\int_D\phi(x)^2\Bigl(F(\nu(t),t)(x)-F(\nu^n_t,t)(x)\Bigr)\,\dif x\Big|
&\leq& L\,\EX^n\int_D |\phi(x)|^2\Big|\int_D w(x,y)\bigl(\nu(t,y)-\nu^n_t(y)\bigr)\dif y\Big|\,\dif x\nonumber\\[1ex]
&\leq& L\,\EX^n\int_D |\phi(x)|^2\, \|w(x,\cdot)\|_{L^2}\,\|\nu(t)-\nu^n_t\|_{L^2}\,\dif x\nonumber\\[1ex]
&\leq& L\,\|\phi\|_{L^4}^2\|w\|_{L^2\times L^2}\,\EX^n\|\nu(t)-\nu^n_t\|_{L^2}.\label{quad_var_proof_est2}
\end{eqnarray}

\medskip

($iii$)\quad In order to estimate the next term we use the bound \eqref{components_mean_bound} on $\EX^n\Theta^{k,n}_t$ and thus obtain
\begin{eqnarray}
\lefteqn{\EX^n\sum_{k=1}^P \frac{\Theta^{k,n}_t}{l(k,n)}\Big|\int_{D_{k,n}} \phi(x)^2\,\dif x-\frac{\rho_n}{l(k,n)}\Bigl(\int_{D_{k,n}}\phi(x)\,\dif x\Bigr)^2\Big|}\nonumber\\
&&\phantom{xxxxx}\leq\ \bigl(1+\|f\|_0\bigr) \sum_{k=1}^P |D_{k,n}|\,\Big|\frac{1}{|D_{k,n}|}\int_{D_{k,n}} \phi(x)^2\,\dif x-
\Bigl(\frac{1}{|D_{k,n}|}\int_{D_{k,n}}\phi(x)\,\dif x\Bigr)^2\Big|\nonumber
\\
&&\phantom{xxxxx\leq\ xxxxxxxxxxxxxxxxx}\mbox{} +\bigl(1+\|f\|_0\bigr)\sum_{k=1}^P |D_{k,n}|\,\Big|1-\frac{\rho_n\,|D_{k,n}|^2}{l(k,n)\,|D_{k,n}|}\Big|
\,\Bigl(\frac{1}{|D_{k,n}|}\int_{D_{k,n}}\phi(x)\,\dif x\Bigr)^2\nonumber\\
&&\phantom{xxxxx}\leq\ \bigl(1+\|f\|_0\bigr)\sum_{k=1}^P \int_{D_{k,n}}\Bigl(\phi(x)-\frac{1}{|D_{k,n}|}\int_{D_{k,n}}\phi(y)\,\dif y\Bigr)^2\dif x\nonumber\\
&&\phantom{xxxxx\leq\ xxxxxxxxxxxxxxxxx}\mbox{}+\bigl(1+\|f\|_0\bigr)\bigg|1-\frac{v_-(n)}{v_+(n)}\frac{\ell_-(n)}{\ell_+(n)}\bigg|\sum_{k=1}^P |D_{k,n}|\,\Bigl(\frac{1}{|D_{k,n}|}\int_{D_{k,n}}\phi(x)\,\dif x\Bigr)^2\,.\nonumber
\end{eqnarray}
Then the estimate is completed applying the Poincar\'e inequality \eqref{poincare_inequality} to the first term, that is, estimating
\begin{equation*}
\sum_{k=1}^P\int_{D_{k,n}}\Bigl(\phi(x)-\frac{1}{|D_{k,n}|}\int_{D_{k,n}}\phi(y)\,\dif y\Bigr)^2\dif x\ \leq \ \frac{\textnormal{diam} (D_{k,n})^2}{\pi^2}\,\|\nabla\phi\|^2_{L^2}\,,
\end{equation*} 
and the observation that the second term is proportional to $\|\nbar\phi^n\|_{L^2(D)}^2$ which is the piecewise constant approximation to $\phi$ based on the partition $\mathcal{D}_n$, see \eqref{overall_proof_pwc_approx}. 
Therefore we overall obtain an upper bound for the difference constituted by the terms ($iii$) in \eqref{the_basic_long_estimate_quad_var_proof} by
\begin{equation}\label{quad_var_proof_est3}
\EX^n\Big|\int_D \phi(x)^2 \nu^n_t(x)\,\dif x - \rho_n\sum_{k=1}^P\frac{\Theta^{k,n}_t}{l(k,n)^2}\,\langle\mathbb{I}_{D_{k,n}},\phi\rangle_{H^{\alpha}}^2\Big| \ \leq \ \delta_+(n)^2\,\frac{1+\|f\|_0}{\pi^2}\,\|\phi\|^2_{H^1}+\bigl(1+\|f\|_0\bigr)\,R(n)\,.
\end{equation}
In the last term
\begin{equation*}
R(n):=\bigg|1-\frac{v_-(n)}{v_+(n)}\frac{\ell_-(n)}{\ell_+(n)}\bigg|\,\|\nbar\phi^n\|_{L^2}^2
\end{equation*}
converges to zero for $n\to\infty$ by assumption \eqref{assumption_on_volumes} and as the sequence $\|\nbar\phi^n\|_{L^2}$ is bounded as it converges to $\|\phi\|_{L^2}$ for $n\to\infty$.
\medskip

($iv$)\quad Finally we consider the difference 
\begin{eqnarray*}
\lefteqn{\EX^n\Big|\int_D \phi(x)^2F(\nu^n_t,t)(x)\,\dif x - \rho_n\sum_{k=1}^P \frac{l(k,n)}{l(k,n)^2}\nbar f_{k,n}(Y^n_t)\,\langle\mathbb{I}_{D_{k,n}},\phi\rangle_{H^{\alpha}}^2
\Big|}\\
&&\phantom{xxxxxxxxxxxxxxxxxxxxxxx}\leq \EX^n\sum_{k=1}^P \Big|\int_{D_{k,n}}\phi(x)^2F(\nu^n_t,t)(x)\,\dif x - \frac{\rho_n}{l(k,n)}\,\nbar f_{k,n}(Y^n_t)\,\langle\phi,\mathbb{I}_{D_{k,n}}\rangle_{H^{\alpha}}^2\Big|\,.
\end{eqnarray*}
We continue estimating the difference in each summand in the final right hand side and obtain using the triangle inequality for the term inside the expectation
\begin{eqnarray*}
\lefteqn{\EX^n\sum_{k=1}^P\Big|\int_{D_{k,n}}\phi(x)^2F(\nu^n_t,t)(x)\,\dif x - \frac{\rho_n}{l(k,n)}\,\nbar f_{k,n}(Y^n_t)\,\langle\mathbb{I}_{D_{k,n}},\phi\rangle_{H^{\alpha}}^2\Big|}\\
&\leq&\underbrace{\EX^n\sum_{k=1}^P\Big|\int_{D_{k,n}}\phi(x)^2\Bigl(F(\nu^n_t,t)(x)-\nbar f_{k,n}(Y^n_t)\Bigr)\,\dif x\Big|}_{(\ast)}+ \underbrace{\EX^n\sum_{k=1}^P\Big|\nbar f_{k,n}(Y^n_t)\Bigr(\int_{D_{k,n}}\phi(x)^2\,\dif x - \frac{\rho_n}{l(k,n)}\,\langle\mathbb{I}_{D_{k,n}},\phi\rangle_{H^{\alpha}}^2\Bigr)\Big|}_{(\ast\ast)}
\end{eqnarray*}

We start with the first term and observe that it possesses the same structure as the term estimated in part (c) of the proof of Theorem \ref{lln} with the only difference that here the function $\phi$ in the integrand is squared. Therefore we obtain the estimate, cf.~\eqref{final_estimate_in_finite_var_proof},
\begin{equation*}
(\ast)\ \leq \ \delta_+(n)\,\frac{L\,\|\phi\|_{L^4}^2}{\tau\pi}\,\Bigl(\sqrt{|D|}\bigl(1+\|f\|_0\bigr)\|\nabla_{\!x}w\|_{L^2\times L^2}+\|\nabla_{\!x} I(t)\|_{L^2}\Bigr)\,.
\end{equation*}

Next, we estimate the second term. Note that $\nbar f_{k,n}$ is bounded by $\|f\|_0$ and thus the remaining term is just as in part ($iii$) of the proof. Hence we obtain the estimate, cf.~\eqref{quad_var_proof_est3},
\begin{equation*}
(\ast\ast)\ \leq\ \delta_+(n)^2\,\frac{\|f\|_0}{\pi^2}\,\|\phi\|_{H^1}^2+ \|f\|_0\,\bigg|1-\frac{v_-(n)}{v_+(n)}\frac{\ell_-(n)}{\ell_+(n)}\bigg|\,\|\nbar\phi^n\|_{L^2}^2\,.
\end{equation*}

Therefore, we overall obtain an upper bound for the difference generated by the terms ($iv$) by
\begin{eqnarray}\label{quad_var_proof_est4}
\lefteqn{\EX^n\Big|\int_D \phi(x)^2F(\nu^n_t,t)(x)\,\mathbb{I}_{[\nu^n_t(x)<1]}\,\dif x - \rho_n\sum_{k=1}^P\frac{1}{l(k,n)}\nbar f_{k,n}(Y^n_t)\,\mathbb{I}_{[\Theta^{k,n}_t<l(k,n)]}\,\langle\mathbb{I}_{D_{k,n}},\phi\rangle_{H^{\alpha}}^2
\Big|}\nonumber\\[1ex]
&&\leq\ \delta_+(n)\,\frac{L\,\|\phi\|_{L^4}^2}{\pi}\,\Bigl(\sqrt{|D|}\bigl(1+\|f\|_0\bigr)\|\nabla_{\!x}w\|_{L^2\times L^2}+\|\nabla_{\!x} I(t)\|_{L^2}\Bigr)+\delta_+(n)^2\,\frac{\|f\|_0}{\pi^2}\,\|\phi\|_{H^1}^2+\|f\|_0\,R(n),\nonumber\\
\end{eqnarray}
where the term $R(n)$ is as in \eqref{quad_var_proof_est3}.\medskip

($v$)\quad To complete the proof we combine the estimates \eqref{quad_var_proof_est1}--\eqref{quad_var_proof_est4} to obtain
\begin{eqnarray*}
\lefteqn{\EX^n\big|\langle G(\nu(t),t)\phi,\phi\rangle_{H^\alpha}-\rho_n\langle G^n(t)\phi,\phi\rangle_{H^\alpha}\big|}\nonumber\\[2ex]
&\leq& \bigl(1+L\,\|w\|_{L^2\times L^2}\bigr)\|\phi\|_{L^4}^2\,\EX^n\|\nu(t)-\nu^n_t\|_{L^2}+ \bigl(1+2\|f\|_0\bigr)\,R(n)\\[2ex]
&&\mbox{} + \delta_+(n)\,\frac{L\,\|\phi\|_{L^4}^2}{\pi}\,\Bigl(\sqrt{|D|}\bigl(1+\|f\|_0\bigr)\|\nabla_{\!x}w\|_{L^2\times L^2)}+\|\nabla_{\!x} I(t)\|_{L^2}\Bigr)  + \delta_+(n)^2\,\frac{1+2\|f\|_0}{\pi^2}\,\|\phi\|_{H^1}^2\,.
\end{eqnarray*}
Integrating over $(0,T)$ we obtain with a suitable constant $C_\phi>0$ independent of $n$ and $T$ the estimate
\begin{eqnarray*}
\lefteqn{\int_0^T\EX^n\big|\langle G(\nu(t),t)\phi,\phi\rangle_{H^\alpha}-\rho_n\langle G^n(t)\phi,\phi\rangle_{H^\alpha}\big|\,\dif t}\\[2ex]
&&\phantom{xxxxx}\leq\ C_\phi\,\Bigl(\EX^n\|\nu(t)-\nu^n_t\|_{L^1((0,T),L^2)} + T\,R(n) + T\,\delta(n)\bigl(1+\|\nabla_{\!x} I(t)\|_{L^1((0,T),L^2)}\bigr)+ \delta(n)^2\Bigr)\,.
\end{eqnarray*}
The constant $C_\phi$ depends on the norm of $\phi$ in the spaces $H^1$ and $L^4$ where the latter can be estimated in terms of the norm in the Sobolev space $H^\alpha$ due to the embedding $H^\alpha\hookrightarrow L^4$, i.e., $C_\phi$ is finite and depends only on $\phi\in H^\alpha$. Finally, each term in the right hand side converges to zero for $n\to\infty$ and hence condition \eqref{clt_martingale_cond_2} follows. The proof of Theorem \ref{martingale_clt} is completed.

\medskip

\textbf{Acknowledgements:} The authors thank J.~Touboul for directing our attention also towards the infinite-time convergence in Theorem \ref{theorem_infinite_time_lln}.

\begin{appendix}

\section{Well-posedness of the Wilson-Cowan equation}\label{section_well_posedness}

This section provides a concise exposition, based on classical existence theory, of the well-posedness of the Wilson-Cowan equation \eqref{Wilson_Cowan} and the boundedness and regularity results for its solution as referred to in Section \ref{section_det_limit}. We understand equation \eqref{Wilson_Cowan} as an $L^2(D)$--valued integral equation, i.e.,
\begin{equation}\label{abstract_Wilson_Cowan}
\nu(t)\,=\, \nu_0 +\frac{1}{\tau} \int_0^t\Bigl(-\nu(s)+F(\nu(s),s)\Bigr)\,\dif s\qquad t\geq 0,\, \nu_0\in L^2(D),
\end{equation}
where the integral is a Bochner integral and $F$ is the Nemytzkii operator acting on $L^2(D)$ defined by
\begin{equation*}
F(g,t)(x)= f\Bigl(\int_D w(x,y)g(y)\,\dif y+I(t,x)\Bigr)\qquad\forall\,g\in L^2(D)\,.
\end{equation*}
As in  Section \ref{section_det_limit} we assume that $f:\rr\to\rr_+$ is Lipschitz continuous, $w\in L^2(D\times D)$ and $I\in C(\rr_+,L^2(D))$, which implies that $F$ is continuous in $t$. Furthermore, it was shown in Section \ref{section_proof_lln} that under these assumptions $F(g,t)$ is Lipschitz continuous in the argument $g$ with Lipschitz constant independent of $t\geq 0$. Thus the integrand in \eqref{abstract_Wilson_Cowan} is Lipschitz continuous with respect to the $L^2(D)$--valued argument for all $t\geq 0$ and, moreover, uniformly continuous in $g$ with respect to $t$. It follows that the integrand in \eqref{abstract_Wilson_Cowan}, that is, the map $(g,t)\to -g+F(g,t)$, is jointly continuous on $\rr_+\times L^2(D)$. Then \cite[Thm.~5.1.1]{Corduneanu} implies that there exists a unique, strongly continuous, global solution to \eqref{abstract_Wilson_Cowan} for every initial condition $\nu_0\in L^2(D)$. By definition this solution is absolutely continuous and, as $F$ is jointly continuous, the derivative of the solution is continuous and exists everywhere. Thus, we conclude that there exists a unique continuously differentiable solution, i.e., $\nu\in C^1(\rr_+,L^2(D))$. 

Next, we recall an `explicit' representation of the solution is the variation of constants formula \eqref{det_sol_var_of_constants} which we already stated in Section \ref{section_proof_infinite_time_lln}. We have that the solution of the Wilson-Cowan equation satisfies the integral equation
\begin{equation*}
\nu(t)=\nu_0+\int_0^t A\nu(t)+ F(\nu(t),t)/\tau\,\dif t,
\end{equation*}
where $A$ is the linear operator in $L^2(D)$ mapping $g$ to $-g/\tau$. Thus, the solution $\nu$ satisfies
\begin{equation*}
\nu(t)=\e^{tA}\nu_0+\frac{1}{\tau}\int_0^t \e^{(t-s)A}F(s,\nu(s))\,\dif s\qquad\forall \,t\geq 0. 
\end{equation*}
In the present setting the application of the linear operator $\e^{tA}$ corresponds to the scalar multiplication with $\e^{-t/\tau}$ as $A=-\tfrac{1}{\tau}\textnormal{Id}_{L^2}$ and thus
\begin{equation*}
\nu(t)=\e^{-t/\tau}\nu_0+\frac{1}{\tau}\int_0^t \e^{-(t-s)/\tau}F(s,\nu(s))\,\dif s\qquad\forall \,t\geq 0. 
\end{equation*}\medskip

We next discuss the results stated in Section \ref{section_det_limit} on the higher spatial regularity of solutions to \eqref{abstract_Wilson_Cowan}. Then a pointwise bound on $\nu(t)\in L^2(D)$, i.e., a constant $C$ such that $|\nu(t,x)|\leq C$ for almost all $x\in D$ and all $t\geq 0$, are then easily obtained by an approximation argument, that is, approximating the less regular solution by solutions of higher regularity. It is possible to prove the pointwise bounds directly, see e.g., \cite{RiedlerPhD} for such an argumentation in a similar setting. However, it is easier and more illustrative to use available results for solutions of higher spatial regularity which are usually arising as the deterministic solution of \eqref{abstract_Wilson_Cowan} one is interested in. E.g., the authors in \cite{FaugerasVeltz} argue that from an application point of view it is reasonable to consider at least continuous solutions. 
In particular, the authors in \cite{FaugerasVeltz} present a detailed existence and uniqueness result for the activity based Amari mean field equation and state that an analogous result hold for the Wilson-Cowan equation \eqref{abstract_Wilson_Cowan} for spatial dimensions $d\leq 3$ which covers all physical relevant domains. Concerning the spatial regularity they consider the space $H^\alpha(D)$, where $\alpha$ is set to be the smallest integer such that $\alpha>d/2$. 
The significance of the choice of $\alpha>d/2$ is -- as so often in this study -- that this implies the embedding of the space $H^\alpha(D)$ into $C(\nbar D)$.  Furthermore we then even obtain that $C([0,T],H^\alpha(D))\subset C([0,T]\times \nbar D)$, i.e., the solution $\nu$ is jointly continuous. 

Therefore we have the subsequent theorem which is sufficient for the set-up in this study. However we note that existence and uniqueness of solutions of the Amari equation were considered under less strict regularity assumptions on the coefficients in \cite{Potthast} and we conjecture that these are also valid for the Wilson-Cowan equation.

\begin{theorem}\label{thm_veltz_faugeras}\textnormal{\cite[Sec.~2]{FaugerasVeltz}} The domain $D$ is bounded and satisfies the strong local Lipschitz property. We assume that $w\in H^{\alpha}(D\times D)$, that $f\in C^\alpha(D)$ with all derivatives bounded, and that $I\in C(\rr_+,H^\alpha(D))$. Then there exists a unique global solution $\nu\in C([0,T],H^\alpha(D))$ for every $T>0$ and every initial condition $\nu_0\in H^\alpha(D)$ to \eqref{abstract_Wilson_Cowan} which depends continuously on the initial condition and is continuously differentiable. Moreover the solution is globally bounded in $H^\alpha(D)$ if the externally applied current $I$ is globally bounded.
\end{theorem}

\begin{remark} In the work \cite{FaugerasVeltz} the authors assume for the domain only the cone property which is implied by the strong local Lipschitz property, see \textnormal{\cite[p.~84]{Adams}}. The latter is the necessary boundary regularity for the present study, cf.~footnote \ref{maurin_footnote} on p.~\pageref{maurin_footnote}. Furthermore, in the reference \cite{FaugerasVeltz} it is also assumed that the gain function  $f$ is infinitely often differentiable with bounded derivatives, but it is surely sufficient for $f$ being $\alpha$-times continuously differentiable.
\end{remark}

Finally, it remains to show the pointwise bound $\nu(t,x)\in (0,\|f\|_0)$ if the initial condition satisfies $\nu_0(x)\in[0,\|f\|_0]$ proposed in Section \ref{section_det_limit}. Under Theorem \ref{thm_veltz_faugeras} the solution $\nu(t,x)$ to \eqref{abstract_Wilson_Cowan} is jointly continuous and therefore the Wilson-Cowan equation holds pointwise in $x$ everywhere and for all $t\geq 0$. Furthermore, $t\mapsto \nu(t,x)$ is continuously differentiable for every fixed $x\in D$ and it is immediate that the bounds are satisfied due to the fact that the derivative of the solution approaching $0$ or $\|f\|_0$ becomes positive or negative, respectively. 
Now, using an approximation result of smooth solutions converging to the $L^2(D)$ solution we obtain that even in this less regular case the pointwise bounds hold almost everywhere.

\section{Comparisons of moment equations}\label{app_moment_equations}

In this section we discuss the moment equations for the $L^2(D)$--valued jump Markov processes $\nu^n_t=\nu^n(\Theta^n_t)$. These can be derived from the corresponding moment equations of the jump Markov process $(\Theta^n_t)_{t\geq 0}$ taking values in $\mathbb{N}^P$. This process is analogous in structure to the usual model used in chemical reaction kinetics, cf., e.g., \cite{Melykuti}. Thus we can use the formulae derived in this reference to obtain, e.g., for the mean the system of differential equations
\begin{equation}\label{mean_eq_MC}
\frac{\dif}{\dif t}\EX^n\Theta^n_t\ =\ -\frac{1}{\tau}\EX^n\Theta^n_t+\frac{1}{\tau}\sum_{k=1}^P l(k,n)\, \EX^n f\Bigl(\sum_{j=1}^P\nbar W_{kj}^n\Theta^{k,n}_t\Bigr)\,e_k\,.
\end{equation}
Furthermore it is straightforward to state a system for the second moments, however, we are not so much interested in the moments of the Markov chain model but those of the $L^2(D)$--valued processes $(\nu^n_t)_{t\geq 0}$ which we can compare to the Langevin approximation. As $\nu^n$ is a linear mapping from $\rr^P$ into $L^2(D)$, it holds that $\nu^n(\EX^n\Theta^n_t)=\EX^n \nu^n(\Theta^n_t)$ and $\nu^n\bigl(\tfrac{\dif}{\dif t}\EX^n\Theta^n_t\bigr)= \tfrac{\dif}{\dif t}\EX^n\nu^n(\Theta^n_t)$, and thus
\begin{equation}
\frac{\dif}{\dif t} \EX^n\nu^n_t = -\frac{1}{\tau}\,\EX^n\nu^n_t + \frac{1}{\tau}\,\EX^n \nbar F^n(\nu^n_t,t)\,.
\end{equation} 

For the second moments of the $L^2(D)$--valued process we obtain for all $\phi\in L^2(D)$
\begin{eqnarray}
\frac{\dif}{\dif t}\EX^n(\phi,\nu^n_t)_{L^2}^2&=& \EX^n\Bigl(\sum_{k=1}^P\frac{\Theta^{k,n}_t}{l(k,n)}\,\int_{D_k} \phi(x)\,\dif x\Bigr)^2\nonumber\\
&=& 2\,\frac{1}{\tau}\,\EX^n\bigl[(\phi,\nu^n_t)_{L^2}\bigl(\phi,-\nu^n_t+\nbar F^n(\nu^n_t,t)\bigr)_{L^2}\bigr]\nonumber\\
&&\mbox{} + \frac{1}{\tau}\,\EX^n\biggl[\sum_{k=1}^P \frac{1}{l(k,n)^2}\Bigl(\Theta^{k,n}_t+l(k,n)\,\nbar f_{k,n}(\nu^n_t,t)\Bigr)\Bigl(\int_{D_{k,n}}\phi(s)\,\dif x\Bigr)^2\Biggr]\nonumber\\[2ex]
&=& \frac{2}{\tau}\,\EX^n\bigl[(\phi,\nu^n_t)_{L^2}\bigl(\phi,-\nu^n_t+\nbar F^n(\nu^n_t,t)\bigr)_{L^2}\bigr] + \EX^n\bigl(G^n(\Theta^n_t,t)\,\phi,\phi\bigr)_{L^2},\label{sec_mom_eq_MC}
\end{eqnarray}
where the bilinear form $\bigl(G^n(\Theta^n_t,t)\,\phi,\phi\bigr)_{L^2}$ is as defined in \eqref{jump_quadratic_var}.\medskip

Next, we state the moment equations for the stochastic partial differential equations. We assume that the Langevin approximation \eqref{langevin_approximation} possesses a (strong) solution in an appropriate Hilbert space $H$ and employ the It{\^o}-formula \cite[Sec.~4.5]{DaPratoZabczyk} which yields for all $\phi\in H^\ast$
\begin{eqnarray*}
\langle\phi,V_t\rangle_{H} &=& \langle\phi,V_0\rangle_{H} + \eps_n\,\int_0^t \bigl\langle\phi,\sqrt{G(V_s,s)}\,\dif W_s\Bigr\rangle_{H}+\int_0^t \bigl\langle\phi,\tfrac{1}{\tau}\,V_s+\tfrac{1}{\tau}\,F(V_s,s)\bigr\rangle_{H}\,\dif s
\end{eqnarray*}
and
\begin{eqnarray*}
\langle\phi,V_t\rangle_{H}^2 &=& \langle\phi,V_0\rangle_{H} + \eps_n\,\int_0^t \bigl\langle 2\langle\phi,V_s\rangle_{H}\,\phi,\sqrt{G(V_s,s)}\,\dif W_s\bigr\rangle_{H}\\[1ex]
&&\mbox{}+2\int_0^t \langle\phi,V_s\rangle_{H}\,\bigl\langle\phi,-\tfrac{1}{\tau}\,V_s+\tfrac{1}{\tau}\,F(V_s,s)\bigr\rangle_{H}\,\dif s+\eps_n^2\,\int_0^t \bigl\langle\phi,G(V_s,s)\phi\bigr\rangle_{H}\,\dif s\,.
\end{eqnarray*}
Next, we take the expectation both sides of these identities and differentiate with respect to $t$ resulting for the first moment in the differential equations 
\begin{equation*}
\frac{\dif}{\dif t}\,\EX\langle\phi,V_t\rangle_{H} = \bigl\langle\phi,\EX\bigl[-\tfrac{1}{\tau}\,V_t+\tfrac{1}{\tau}\,F(V_t,t)\bigr]\bigr\rangle_{H}
\end{equation*}
which is equivalent to the abstract evolution equation in $H$ given by
\begin{equation}\label{mean_eq_langevin}
\frac{\dif}{\dif t}\,\EX V_t = -\frac{1}{\tau}\,\EX V_t+\frac{1}{\tau}\,\EX F(V_t,t)\,.
\end{equation}
And for the second moment we obtain the differential equation
\begin{eqnarray}\label{sec_mom_eq_langevin}
\frac{\dif}{\dif t}\,\EX\langle\phi,V_t\rangle_{H}^2&=&\frac{2}{\tau}\,\EX\bigl[\langle\phi,V_t\rangle_{H}\,\bigl\langle\phi,-V_t+F(V_t,t)\bigr\rangle_{H}\bigr]+\eps_n^2\,\EX\bigl[\bigl\langle G(V_t,t)\,\phi,\phi\bigr\rangle_{H}\bigr]\,.
\end{eqnarray}
Further, the linear noise approximation \eqref{linear_noise_approximation} satisfies the equations
\begin{equation}\label{mean_eq_lin}
\frac{\dif}{\dif t}\,\EX U_t = -\frac{1}{\tau}\,\EX U_t+\frac{1}{\tau}\,\EX F(U_t,t)\,.
\end{equation}
and
\begin{eqnarray}\label{sec_mom_eq_lin}
\frac{\dif}{\dif t}\,\EX\langle\phi,U_t\rangle_{H}^2&=&\frac{2}{\tau}\,\EX\bigl[\langle\phi,U_t\rangle_{H}\,\bigl\langle\phi,-U_t+F(U_t,t)\bigr\rangle_{H}\bigr]+\eps_n^2\,\bigl\langle G(\nu(t),t)\,\phi,\phi\bigr\rangle_{H}\,.
\end{eqnarray}
Finally, we note that exactly the same moment equations hold for the variants of the linear noise and Langevin approximation using a $Q$--Wiener process and an appropriate diffusion coefficient, cf.~Remark \ref{remark_non-uniqueness_spde_approx}.\medskip

A comparison of the moment equations \eqref{mean_eq_MC}, \eqref{mean_eq_langevin}, \eqref{mean_eq_lin} for the mean and \eqref{sec_mom_eq_MC}, \eqref{sec_mom_eq_langevin}, \eqref{sec_mom_eq_lin} for the second moments show that they are similar in structure but do not coincide. This is analogous to the properties of the moment equations in finite dimension and as in finite dimensions there is one exception, which is the case of first order transitions:
If $F$ were affine\footnote{In the case of an affine function $F(v,t)=f_1(t)\cdot v+f_2(t)$ the mapping $f_1$ is a linear form on $H$ which is interchangeable with the expectation operator. For example, in the simplest case the application of $f_1$ to $v\in H$ is just a multiplication by a scalar.} in $v$, i.e., $F(v,t)=f_1(t)\cdot v+f_2(t)$, then we obtain that the first moment equations \eqref{mean_eq_langevin} and \eqref{mean_eq_lin} of the Langevin and linear noise approximation, respectively, reduce to the Wilson-Cowan equation with $\nu(t)=\EX V_t=\EX U_t$. 
Furthermore, if $F$ is affine, this implies that also $G$ is affine in $V_t$ and thus
\begin{equation}
\langle\phi, G(V_t,t)\phi\rangle_H = \frac{1}{\tau} \langle\phi,V_t\cdot\phi\rangle_H+\langle\phi,f_1(t)\cdot V_t\cdot \phi \rangle_H + \langle\phi, f_2(t)\cdot\phi\rangle_H\,.
\end{equation}
Taking the expectation on both sides and assuming interchangeability of the expectation with the application of all the linear forms (think of the duality pairing as the inner product in $L^2(D)$) we obtain
\begin{equation}
\EX\langle G(V_t,t)\,\phi,\phi\rangle_H = \frac{1}{\tau} \langle\phi,\EX[V_t]\cdot\phi\rangle_H+\langle\phi,f_1(t)\cdot \EX[V_t]\cdot \phi \rangle_H + \langle\phi, f_2(t)\cdot\phi\rangle_H = \langle G\bigl(\EX[V_t],t\bigr)\,\phi,\phi\rangle_H\,.
\end{equation}
As $\EX V_t =\EX U_t=\nu(t)$ we obtain that the second moment equation for the Langevin approximation and the linear noise approximation coincide. Moreover, they are closed (for each $\phi$), i.e., the system depends only on $\EX V_t$ and $\EX\langle\phi,V_t\rangle_H^2$. Again, this corresponds to the well-known case from finite-dimensional chemical reaction kinetics.

\medskip

Finally, if $F$ is affine also the connection of the moment equations for the resulting Markov chain models is interesting. On the one hand the equation for the mean coincides with the Wilson-Cowan equation where the gain function in its right hand side is given by $\nbar F^n$. As $\nbar F^n$ is essentially a piecewise constant approximation to $F$ the resulting equations for the mean correspond to a spatial discretisation of the Wilson-Cowan equation, cf.~the continuum limit in the derivation of the mean field equation in \cite{Bressloff1}.

\end{appendix}

\bibliographystyle{plain}
\bibliography{bibliography}

\end{document}